\documentclass[12pt]{amsart}
\usepackage{fancyhdr}
\usepackage{amsmath}
\usepackage{amsthm}
\usepackage{latexsym}
\usepackage{amssymb}
\DeclareMathOperator{\rank}{rank}

\DeclareMathOperator{\ad}{ad}
\DeclareMathOperator{\sgn}{sign}

\newtheorem{theorem}{Theorem}
\newtheorem{proposition}{Proposition}
\newtheorem{lemma}{Lemma}

\newtheorem{definition}{Definition}

\newcounter{obsctr}
\newenvironment{observation}{\stepcounter{obsctr}%
\vskip 4pt \noindent {\bf \arabic{obsctr}.}\hskip 5pt }{\vskip 4pt}
\newtheorem{remark}{Remark}
\renewcommand{\thetheorem}{\thesection.\arabic{theorem}}
\renewcommand{\theproposition}{\thesection.\arabic{proposition}}
\renewcommand{\thelemma}{\thesection.\arabic{lemma}}
\renewcommand{\thedefinition}{\thesection.\arabic{definition}}
\renewcommand{\thecorollary}{\thesection.\arabic{corollary}}
\renewcommand{\theequation}{\thesection.\arabic{equation}}
\begin{document}
\baselineskip 16pt
\def\R {{\mathbb{R}}}
\def\N {{\mathbb{N}}}
\def\C {{\mathbb{C}}}
\def\Z {{\mathbb{Z}}}
\def\phi{\varphi}
\def\epsilon{\varepsilon}
\title{A Class of Sums of Squares with a Given Poisson-Treves
  Stratification}    
\author{Antonio Bove}
\address{Dipartimento di Matematica, Universit\`{a} di Bologna,
Piazza di Porta S. Donato 5, 40127 Bologna, Italia, and Istituto Nazionale
di Fisica Nucleare, Sezione di Bologna}
\email{bove@dm.unibo.it}
\author{David Tartakoff}
\address{Department of Mathematics, University
of Illinois at Chicago, m/c 247, 851 S.
Morgan St., Chicago IL  60607}
\email{dst@uic.edu}
\date{\today}
\begin{abstract}
We study a class of sum of squares exhibiting the same Poisson-Treves
stratification as the Oleinik-Radkevi\v c\ operator. We find three
types of operators having distinct microlocal structures. For one of
these we prove a Gevrey hypoellipticity theorem analogous to our
recent result for the corresponding Oleinik-Radkevi\v c\ operator.
\end{abstract}
\maketitle
\fancyhead{}
\fancyhead[RE]{\small A. Bove and D.S. Tartakoff}
\fancyhead[LE]{\thepage}
\fancyhead[LO]{\small A Class of Sums of Squares}
\fancyhead[RO]{\thepage}
\fancyfoot[CO,CE]{}
\renewcommand{\headrulewidth}{0pt}
\renewcommand{\footrulewidth}{0pt}
\section{Introduction}
\renewcommand{\theequation}{\thesection.\arabic{equation}}
\setcounter{equation}{0}
\setcounter{theorem}{0}
\setcounter{proposition}{0}  
\setcounter{lemma}{0}
\setcounter{corollary}{0} 
\setcounter{definition}{0}

The problem of analytic hypoellipticity for second order
operators which are sums of squares of vector fields with
analytic coefficients has been widely studied and has been
around since the paper of L.~H\"or\-man\-der \cite{H1967} on $
C^{\infty} $ hypoellipticity for this type of operator. In
particular D.~S.~Tartakoff
\cite{Tartakoff1978}, \cite{Tartakoff1980}, for second order, and F.~Treves
\cite{Treves1978} for general order, gave general analytic
hypoellipticity theorems for the case when the
characteristic manifold is symplectic and the operator degenerates on
it to an exact order. In the case of non exact and higher order
degeneracy O.~Oleinik \cite{O 1973} and O.~Oleinik and R.~Radkevi\v c
\cite{OR1974} (see also Christ \cite{Christ1996}) showed that in
general one
cannot have analytic hypoellipticity (see also the paper \cite{BT} by
the present authors for a precise and optimal partial regularity result in
the case of the operator studied by Oleinik and Radkevi\v c), but only 
certain degrees of Gevrey hypoellipticity.

Into this scenario there appeared in 1999 the well known paper by F.~Treves
\cite{trevesprepr} introducing the notion of Poisson stratification for
a set of vector fields satisfying H\"ormander's Lie algebra condition
and having analytic coefficients. Basically, crudely simplifying
Treves' setting, the conjecture states that an operator which is a sum of
squares of vector fields is analytic hypoelliptic if and only if every
layer in its Poisson stratification is symplectic. To our knowledge the
conjecture has been neither proved nor disproved up to now.

In this paper, the first of a series, we study an operator which is the
sum of the squares of three vector fields with analytic coefficients
in three variables. For such an operator we assume that its
Poisson-Treves stratification is given in such a way that its
H\"ormander numbers are the same as those of the Oleinik-Radkevi\v c
model operator. Here by H\"ormander numbers we mean both the number
and the relative codimensions of the stratification's layers. Our
purpose is to classify such kind of operators and obtain their Gevrey
(or possibly analytic) hypoellipticity threshold.

In the present paper we classify the operators  having the required
properties and, for one of the classes, we obtain the same Gevrey
hypoellipticity threshold as that of the Oleinik-Radkevi\v c model
operator. We are unable to deduce these (optimal)
results for every class of operators sharing the Poisson-Treves
stratification with the Oleinik-Radkevi\v c model, but we shall come
back to this subject in a forthcoming paper \cite{BTforth}.

Before stating our assumptions precisely, we want to make a
couple of remarks. 

\begin{observation} Our vector fields are linearly independent outside of the
characteristic manifold. This essentially implies that the
characteristic manifold is cylindrical with respect to a two
dimensional subspace of the fibers of the cotangent bundle, or in
other words, it is the zero set of one covariable and one function of
the variables in the base. This restriction eliminates cases where two of
the vector fields can become colinear outside of the characteristic
set. On the other hand, many results are known for the case
of the sums of two squares.
\end{observation}

\begin{observation} The Oleinik-Radkevi\v c model operator has a
codimension 2 symplectic characteristic manifold. In three dimensions one
might consider also cases where the characteristic manifold is symplectic
and of dimension 4 or has symplectic layers of codimension 2 and
symplectic layers of codimension 4. Even for the sums of two squares,
though, this situation faces difficulties of the same kind as those appearing
in Christ's example \cite{Christ1998}.
\end{observation}

In the first part of
the paper we deduce some standard forms (cf. Theorem \ref{th1}) below)
that can be useful in proving  {\em a priori} estimates. Then we proceed
to prove a Gevrey hypoellipticity threshold for one of these standard
forms (cf. Theorem \ref{theorem2})

Essentially the
operators verifying our assumptions fall into three classes, depending on
how the vector fields vanish on the characteristic set. For the first
case, called Case 1, we make a finer analysis of the extent to which the
vector fields under consideration are linearly independent outside of the
characteristic manifold. This is accomplished by looking at each of the
``characteristic'' vector fields and computing it on the null
bicharacteristic curve of the only non characteristic vector field.

This operation does not affect the covariables (i.e. affects only the
coefficients of the base), since the null bicharacteristic curve is a
curve in the base variables. Then one focusses on the zero set of the
resulting vector fields. Computing the symbol of one of the
vector fields on the zero set of the other allows us to define a sort of
degeneracy rate which turns out to be useful in the {\em a priori}
estimates. The last section of the present paper is concerned with the
case when the above mentioned degeneracy rate is zero. Then we obtain the
same (optimal) result as for the Oleinik-Radkevi\v c model.

If the degeneracy rate is larger than zero, the estimates are deduced
in a very different way and yield a different result. This is the
subject of a forthcoming paper.

The second and third classes (called Case 2a and Case 2b below)
will also be studied in a subsequent paper.

%
%
\section{Assumptions}
\renewcommand{\theequation}{\thesection.\arabic{equation}}
\setcounter{equation}{0}
\setcounter{theorem}{0}
\setcounter{proposition}{0}  
\setcounter{lemma}{0}
\setcounter{corollary}{0} 
\setcounter{definition}{0}

We now specify the assumptions.
Our operators have the general form 
$$ 
P (x, D) = \sum_{j=1}^{3} X_{j}^{2}(x, D) ,
$$
where $ x \in \R^{3} $ and $ D_{j} = \frac{1}{i}
\frac{\partial}{\partial x_{j}} $, $ j = 1, 2, 3 $. Here $ X_{j} $
denotes a vector field with real analytic coefficients defined in a
neighborhood of the origin in the $ x $ variable.

The following assumptions try to mimic the fact that $ P $ has the
same Poisson-Treves stratification as the operator $ D_{1}^{2} +
x_{1}^{2(p-1)} D_{2}^{2} + x_{1}^{2(q-1)} D_{3}^{2} $, where $ p $ and
$ q $ are integers and $ 1 \leq p \leq q $.
\smallskip

\begin{itemize}
\item[(A1)]{}
The operator $ P $ satisfies the H\"ormander Lie algebra condition and
hence is $ C^{\infty} $ hypoelliptic. As a consequence not all the
vector fields are characteristic (i.e. have vanishing coefficients) on
the characteristic manifold. Hence we may suppose without loss of
generality that 
$$ 
X_{1}(x, D) = D_{1}.
$$
\vskip .5cm
\item[(A2)]{}
We may always assume that the point $ (0; e_{3}) $ is a characteristic
point for $ P $ (using a translation and a rotation if necessary). We
assume then 
that near $ (0; e_{3}), $ the characteristic set of $ P $ is an
analytic symplectic submanifold of
$T^{*}\R^{3}\setminus 0$ of codimension two which we denote by $ \Sigma_{1}.
$  We explicitly note that this is a microlocal
assumption.
\vskip .5cm 
\item[(A3)]{}
Let $ \Omega = U \times \Gamma $ a conic neighborhood of the
point  $ (0, e_{3}) $. And let $ \pi_{1} \colon U \times \Gamma
\rightarrow U $ be the projection onto the space variables. We assume that
the vector fields
$$ 
\left . {X_{j}} \right |_{U \setminus \pi_{1}\Sigma_{1}}
$$
are linearly independent; the above notation means that restricting
the coefficients of the fields $ X_{j} $ to the space projection of $
\Sigma_{1} $ yields linearly independent vectors in $ \R^{3} $.

Note that, because of (A2), the coefficients of the vector fields
depend non trivially on the $ x $ variable. This assumption has strong
implications on the structure of $ \Sigma_{1} $ and, to avoid
technical details at this point, we refer to Section 3.
\vskip .5cm
\item[(A4)]{}
\begin{multline*} 
\qquad\quad\Sigma_{2} = \{ (x, \xi) \in T^{*}\R^{3}\setminus 0 \ | \ (x,
\xi)
\in
\Sigma_{1}, \{ X_{i}, X_{j} \}(x, \xi) = 0, 
\\
\ i,j \in \{1,2,3\} \},
\end{multline*}
and, in general, let $ I = (i_{1}, \ldots , i_{k}) $, $ i_{j} \in \{1,
2, 3 \}$, for $ j = 1, \ldots k $. Writing $ | I | = k, $ we denote
by $ X_{I} $ the iterated Poisson bracket
$$ 
X_{I} = \{ X_{i_{1}} ,\ \{ X_{i_{2}},\ \ldots \ , \{ X_{i_{k-1}},\
X_{i_{k}} \} \ \ldots \} \}
$$
of the vector fields $ X_{j} $, $ j = 1, 2, 3 $; set
\begin{multline*}
\qquad\quad\;\Sigma_{h} = \{ (x, \xi) \in T^{*}{\R}^{3} \setminus 0 \
| \ (x,
\xi) \in \Sigma_{h-1}, \ X_{I}(x, \xi) = 0 \\ 
\text{for every index $ I $ such that $ |I| = h $} \}.
\end{multline*}
Let $ p\leq q $ be two positive integers. Then
we make the following assumptions:
\begin{itemize}
\item[(i)]{}
$ \Sigma_{1} \cap \Omega = \cdots = \Sigma_{p-1} \cap \Omega $.
\item[(ii)]{}
$ \Sigma_{p} \cap \Omega $ is a non-empty analytic submanifold of $
\Sigma_{1} \cap \Omega $ of codimension one.
\item[(iii)]{}
$ \Sigma_{p} \cap \Omega = \Sigma_{p+1} \cap \Omega = \cdots =
\Sigma_{q-1} \cap \Omega $.
\item[(iv)]{}
$ \Sigma_{q} \cap \Omega $ is empty in $ T^{*}{\R}^{3} \setminus 0
$ (i.e. $ \Sigma_{q} \cap \Omega $ is contained in the zero section of
the cotangent bundle over $ \Omega $.)
\end{itemize}
\end{itemize}

\section{Standard Forms: The equations of $\Sigma_1$}
\setcounter{equation}{0}
\setcounter{theorem}{0}
\setcounter{proposition}{0}  
\setcounter{lemma}{0}
\setcounter{corollary}{0} 
\setcounter{definition}{0}

Due to the above assumptions we may suppose that the vector fields have
the following form:
\begin{eqnarray}
\label{3.1}
X_{1}(x, \xi) & = & \xi_{1} \nonumber\\
X_{2}(x, \xi) & = & a_{21}(x) \xi_{1} + a_{22}(x) \xi_{2} + a_{23}(x)
\xi_{3}  \\
X_{3}(x, \xi) & = & a_{31}(x) \xi_{1} + a_{32}(x) \xi_{2} + a_{33}(x)
\xi_{3}. \nonumber
\end{eqnarray}
Hence $ \xi_{1} = 0  $ is one of the two equations defining $ \Sigma_{1}
$; letting 
$$ 
A(x) = \begin{bmatrix}
  a_{22}(x) & a_{23}(x) \\
  a_{32}(x) & a_{33}(x)
\end{bmatrix},
$$
(the $a_{jk}$ being analytic), the other equation is given by
\begin{equation}
\label{3.2}
A(x) \xi' = 0,
\end{equation}
where $ \xi' = (\xi_{2}, \xi_{3}). $

We claim that this can only be the second defining condition of
$\Sigma_1$ if $A(x)\equiv 0$ on $\Sigma_1$ (locally). For suppose
$ (x_{0},
\xi_{0}')
$,
$
\xi_{0}'
\neq 0
$, is such that
\begin{equation}
\label{3.3}
A(x_{0}) \xi_{0}' = 0
\end{equation}
and assume that to the contrary, for $x$ near $x_0$ on $\Sigma_1,$
$$ 
A(x) \neq \begin{bmatrix}
        0 & 0 \\
  0 & 0
\end{bmatrix} \qquad ({\hbox{which we will write as }} A(x)\neq 0).
$$
\noindent 
This implies that at $x_0$ the rank of $ A $ is equal to $ 1 $ since
$ 0
\neq \xi_{0}' \in \ker A(x_{0}). $
It follows that, in a conic neighborhood of $ (x_{0}; 0, \xi_{0}'), $
the characteristic manifold $ \Sigma_{1} $ is defined by
$$ 
\Sigma_{1} = \{ (x, \xi) \ | \ \xi_{1} = 0,
\det A(x) = 0, \xi' \in \ker A(x) \}, 
$$
because we may always assume that $ \rank A(x) \geq 1 $ near $ x_{0}
$.

Since $ A(x_{0}) \neq 0 $, the latter
two equations in the definition of $ \Sigma_{1} $ are certainly
independent (the second of them has non-zero $ \xi' $-gradient, while
the first of them must
have a non-zero $ x $-gradient).
As a consequence one of them must be identically
satisfied in order to accomplish the codimension 2 condition. 
Since $ \rank A(x_{0}) = 1 $, the condition $
\xi' \in \ker A(x) $ cannot be identically satisfied. 
Hence the only possibility is that
$$ \det A(x) \equiv 0 $$ in a full neighborhood of $ x_{0}
$.
However this fact would imply that there exist points $ (x, \xi) $, $
\xi_{1} = 0 $, $ (x, \xi) \notin \Sigma_{1} $, such that the vector
fields $ X_{1} $, $ X_{2} $, $ X_{3} $ are not linearly independent.

Consequently the only possible case left is that $A$ is the zero matrix:
\begin{equation}
\label{3.4}
A(x) = 0, 
\end{equation}
if $ (x, \xi) \in \Sigma_{1} $. This means that
\begin{equation}
\label{3.5}
\Sigma_{1} = \{ (x, \xi) \ | \ \xi_{1} = 0, A(x) = 0 \}.
\end{equation}
Hence the matrix condition $ A(x) = 0 $ must be (locally) equivalent to $
\phi(x) = 0
$, where $ \phi $ is a real analytic scalar function and
such that $ d_{x}\phi(x) \neq 0. $

By Assumption (A2), $ \{ \xi_{1}, \phi(x) \} \neq 0 $ at $
\Sigma_{1} $. Hence by the implicit function theorem the equation $
\phi(x) = 0 $ is equivalent to the equation
\begin{equation}
\label{3.6}
x_{1} - g(x') = 0,
\end{equation}
where $ g $ is a suitable real analytic function, $ x' = (x_{2},
x_{3}) $, and $ g $ is defined locally.
 We conclude then that
\begin{equation}
\label{3.7}
\Sigma_{1} = \{ (x, \xi) \ | \ \xi_{1} = 0,  \; x_{1} - g(x') = 0 \},
\end{equation}
and that
\begin{equation}
\label{3.8}
A(x) = ( x_{1} - g(x')) \tilde{A}(x),
\end{equation}
for a suitable $ 2 \times 2 $ matrix $ \tilde{A} $ with real analytic
entries $ \tilde{a}_{ij} $, $ i,j \in \{2, 3\} $.

Next we perform a change of variables (and hence a canonical
transformation) which is linear in $ \xi $, so that vector fields are
mapped to vector fields in the new coordinates, allowing us to make the
function $ g $ identically zero.

Define:
\begin{equation}
\label{3.9}
\begin{array}{rclcrcl}
y_{1} & = & x_{1} + g(x') & \qquad & \eta_{1} & = & \xi_{1} \\[7pt]
y_{2} & = & x_{2}         & \qquad & \eta_{2} & = & \xi_{2} - \xi_{1}
\frac{\partial g}{ \partial x_{2}} \\[7pt]
y_{3} & = & x_{3}         & \qquad & \eta_{3} & = & \xi_{3} - \xi_{1}
\frac{\partial g}{\partial x_{3}}.
\end{array}
\end{equation}
The three vector fields become:
\vskip.2in
$$ 
\begin{array}{rcl}
X_{1}(y, \eta) & = & \eta_{1} \\[7pt]
X_{2}(y, \eta) & = & \left( a_{21}(y_{1}-g(y'), y') +
  y_{1}\tilde{a}_{22}(y_{1}-g(y'), y') \frac{\partial g(y')}{\partial
    y_{2}} \right .\\[7pt]
               &   & \left .
+ y_{1} \tilde{a}_{23}(y_{1}-g(y'), y') \frac{\partial
    g(y')}{\partial y_{3}} \right) \eta_{1} \\[7pt]
               &   &  + 
y_{1} \left [\tilde{a}_{22}(y_{1}-g(y'), y') \eta_{2} +
  \tilde{a}_{23}(y_{1}-g(y'), y') \eta_{3} \right] \\[7pt]
X_{3}(y, \eta) & = & \left( a_{31}(y_{1}-g(y'), y') +
  y_{1}\tilde{a}_{32}(y_{1}-g(y'), y') \frac{\partial g(y')}{\partial
    y_{2}} \right .\\[7pt]
               &   & \left .
+ y_{1} \tilde{a}_{33}(y_{1}-g(y'), y') \frac{\partial
    g(y')}{\partial y_{3}} \right) \eta_{1} \\[7pt]
               &   &  + 
y_{1} \left [\tilde{a}_{32}(y_{1}-g(y'), y') \eta_{2} +
  \tilde{a}_{33}(y_{1}-g(y'), y') \eta_{3} \right].
\end{array}
$$
The above fields can be rewritten, with obvious notation, in the following
way:
\begin{equation}
\label{3.10}
\begin{array}{rcl}
X_{1}(x, \xi) & = & \xi_{1} \\[7pt]
X_{2}(x, \xi) & = & a_{21}(x) \xi_{1} + x_{1} \left[ a_{22}(x)\xi_{2}
  + a_{23}(x) \xi_{3}\right] \\[7pt]
X_{3}(x, \xi) & = & a_{31}(x) \xi_{1} + x_{1} \left[ a_{32}(x)\xi_{2}
  + a_{33}(x) \xi_{3}\right]
\end{array}
\end{equation}
with suitable real analytic functions $ a_{ij} $ defined in a
neighborhood of the origin.
\smallskip

\section{Standard forms: the equations of $\Sigma_2, \ldots \Sigma_{p-1}$}
\renewcommand{\theequation}{\thesection.\arabic{equation}}
\setcounter{equation}{0}
\setcounter{theorem}{0}
\setcounter{proposition}{0}  
\setcounter{lemma}{0}
\setcounter{corollary}{0} 
\setcounter{definition}{0}

Let us now turn to Assumption (A4) concerning $ \Sigma_{2} $. We
have
$$ 
\{X_{1}(x, \xi) , X_{j}(x, \xi) \} = \frac{\partial}{\partial x_{1}}
X_{j}(x, \xi),
$$
for $ j = 2, 3 $ and the latter quantity is equal to
$$ 
\frac{\partial a_{j1}(x)}{\partial x_{1}} \xi_{1} + 
\left[ a_{j2}(x)\xi_{2} + a_{j3}(x) \xi_{3}\right] + O(|x_{1}|),
$$
for $ j = 2, 3 $, and
\begin{multline*}
\{ X_{2}, X_{3} \}(x, \xi)  \\
 = \{a_{21}\xi_{1} + x_{1} \left[ a_{22}(x)\xi_{2} + a_{23}(x)
  \xi_{3}\right] , a_{31}\xi_{1} + x_{1}  \left[ a_{22}(x)\xi_{2}
  + a_{23}(x) \xi_{3}\right] \},
\end{multline*}
which gives
\begin{multline}
\label{4.1}
\{ X_{2}, X_{3} \}(x, \xi) = a_{21}(x) \{X_{1}, X_{2}\}(x, \xi) -
a_{31}(x) \{X_{1}, X_{3}\}(x, \xi) \\
+ O(|x_{1}| + |\xi_{1}|),
\end{multline}
where $ O(|x_{1}| + |\xi_{1}|) $ stands for a vector field with
principal symbol vanishing on $ \Sigma_{1} $.
Hence we obtain $ \Sigma_{2} = \Sigma_{1} \cap \{ (x, \xi) | \
\{X_{1}, X_{j} \} = 0 \} $, $ j= 2, 3 $.

Let us again denote by $ A(x) $ the $ 2 \times 2 $ matrix
$$ 
A(x) =
\begin{bmatrix}
  a_{22}(x) & a_{23}(x)  \\
  a_{32}(x) & a_{33}(x)
\end{bmatrix};
$$
then Assumption (A4) means that
\begin{equation}
\label{4.2}
A(x) \xi' = 0
\end{equation}
if and only if $ x_{1} = 0 $. This implies that 
\begin{equation}
\label{4.3}
A(x) = x_{1} \tilde{A}(x)
\end{equation}
for a suitable $ 2 \times 2 $ matrix $ \tilde{A} $ with analytic
entries. 

Iterating the above argument we can conclude that the vector fields
can be written in the form
\begin{equation}
\label{4.4}
\begin{array}{rcl}
X_{1}(x, \xi) & = & \xi_{1} \\[7pt]
X_{2}(x, \xi) & = & a_{21}(x) \xi_{1} + x_{1}^{p-1} \left[ a_{22}(x)\xi_{2}
  + a_{23}(x) \xi_{3}\right] \\[7pt]
X_{3}(x, \xi) & = & a_{31}(x) \xi_{1} + x_{1}^{p-1} \left[ a_{32}(x)\xi_{2}
  + a_{33}(x) \xi_{3}\right].
\end{array}
\end{equation}
We summarize what has been proved up to this point in the 
\begin{proposition}
Suppose that (A1)--(A3) and (A4)(i) hold. Then the vector fields $
X_{1} $, $ X_{2} $, $ X_{3} $ can be written, in a suitable system of
local coordinates, in the form (\ref{4.4}).
\end{proposition}

\section{The equation defining $ \Sigma_{p} $ with respect to $
  \Sigma_{1}$} 
\setcounter{equation}{0}
\setcounter{theorem}{0}
\setcounter{proposition}{0}  
\setcounter{lemma}{0}
\setcounter{corollary}{0} 
\setcounter{definition}{0}

Let us denote by $ \phi (x', \xi')  $ a real analytic function defined
on a (conic) neighborhood of $ (0, e_{3}) $ in $ \Sigma_{1} $ and such
that $ d_{(x', \xi')}\phi(0, e_{3}) \neq 0 $ and the equation $ \phi
(x', \xi') = 0 $ is equivalent to $ A(0, x')\xi' = 0 $. 

We have either
\begin{equation}
\label{5.1}
\frac{\partial \phi}{\partial \xi'}(0, e_{3}) \neq 0
\qquad
\text{ (Case I) } \end{equation}
or
\begin{equation}
\label{5.2}
\frac{\partial \phi}{\partial x'}(0,
e_{3}) \neq 0 \qquad\text{ (Case II)}.
\end{equation}
\subsection{Case I}
For the case of non-zero $\xi$ gradient, assume that it is the
$\xi_2$ derivative of $\phi$ that is non-zero at $(0,e_3)$ (we will see below
that the case of a non-zero
$\xi_3$ derivative cannot occur). 
Then we may write
\begin{equation} 
\label{5.3}
\phi (x', \xi') =  (\xi_{2} - \chi(x', \xi_{3}))e (x', \xi'),
\end{equation}
where $ e $ and $ \chi $ are analytic and $ e(0, e_{3}) \neq 0 $
and thus 
$$ \phi (x', \xi') = 0 \Longleftrightarrow  \xi_{2} -
\chi(x', \xi_{3}) = 0 \Longleftrightarrow A(x')\xi' = 0. $$ 

We claim that $\chi(x',\xi')$ has the simpler form
$\tilde\chi(x')\xi_3,$ and to see this let
$ t $ denote a non-zero real number; if
$ A(x') \xi' = 0 $ then obviously $ A(x') t\xi' = 0 $. Thus
$ (x', \xi')  \in \Sigma_{2} \implies  (x', t \xi')
\in \Sigma_{2} $, so that $ t \xi_{2} - \chi(x', t \xi_{3})
= 0$. Since $ \xi_{2} = \chi(x', \xi_{3}) $,
we have
$ 
\chi(x', t \xi_{3}) = t \chi(x', \xi_{3})
$
for every non-zero real number $ t $.
But now $ \xi_{3} \neq 0 $ in a conic neighborhood of $ (0, e_{3}),$ 
so that
$ \chi(x', \xi_{3}) = \xi_{3} \chi(x', 1) = \xi_{3} \tilde{\chi}(x'),$
for a suitable analytic function $ \tilde{\chi} $ of the space
variable only, and so finally we obtain
\begin{equation}
\label{5.4}
A(x') \xi' = 0 \Longleftrightarrow \xi_{2} - \chi(x') \xi_{3} = 0,
\end{equation}
where we have written $ \chi $ again for the function $ \tilde{\chi} $.

The above formula has been derived in the case that it is the $\xi_2$
derivative of $\phi$ that is non-zero at $(0,e_3).$ Now suppose that 
the $\xi_3$ derivative of $\phi$ is non-zero at $(0,e_3)$ instead.
Then arguing as above we find that the equation $ \phi (x', \xi') =
0 $ is equivalent to $ \xi_{3} - \chi(x', \xi_{2}) = 0 $. As
before, let $ t $ be a non-zero real number; now, since if $ (x', \xi') 
$ belongs to $ \Sigma_{p} $ then also $ (x', t \xi')  $ belongs to $
\Sigma_{p} $, keeping in mind that, by assumption, the point $
(0, e_{3}) $ belongs to $ \Sigma_{p} $, we find that $ t = \chi(0, 0)
$ for any $ t \in \R \setminus 0$, which is absurd.

We thus have proved that if (\ref{5.1}) is true then the equation
defining $ \Sigma_{p} $ relatively to $ \Sigma_{1} $ is given by
(\ref{5.4}). 

\subsection{Case II}

We now turn to the case where $ \phi_\xi'(0,
e_{3}) = 0$ but
\begin{equation}
\label{5.5} 
\frac{\partial \phi}{\partial x'}(0, e_{3}) \neq 0 \qquad\qquad
\hbox{(Case II)}
\end{equation}
and we assume here that 
\begin{equation} 
\label{5.6}
\frac{\partial \phi}{\partial x_{2}}(0, e_{3}) \neq 0 \qquad\qquad 
\hbox{(Case II${}_{x_2}$)}.
\end{equation}
The case $\phi_{x_3}(0, e_{3}) \neq 0 $
(Case II${}_{x_3}$) has some obvious but non-trivial differences that
we shall stress later.

Arguing along the same lines as above we obtain that there is a
function $ \chi(x_{3}, \xi') $ such that the equation $ \phi (x', \xi') =
0$ is equivalent to $ x_{2} - \chi(x_{3}, \xi') = 0 $. Here $ \chi $ is
analytic and defined on a conic neighborhood of $ (0, e_{3}) $ in $
\R_{x_{3}} \times (\R^{2}_{\xi'}\setminus 0) $. 
Again we may assume that on that neighborhood $ \xi_{3} $ is not
zero. Moreover if $ t $ denotes a non-zero real number we obtain that
$ \chi(x_{3}, t \xi') = \chi(x_{3}, \xi') $, so that 
\begin{equation}
\label{5.7}
A(x') \xi' = 0 \Longleftrightarrow x_{2} - \chi(x_{3},
\frac{\xi_{2}}{\xi_{3}}) = 0,
\end{equation}
where we have denoted by $ \chi(x_{3}, \sigma) $ the function $
\chi(x_{3},
\sigma, 1)$.

We point out that the function $ \chi $ in (\ref{5.7}) is an analytic
function defined in a neighborhood of the origin in $
\R_{x_{3}} \times \R_{\sigma} $.

From (\ref{5.7}) we obtain that there is a positive integer $ k $ such
that 
\begin{equation}
\label{5.8}
A(x') \xi' = \left( x_{2} - \chi(x_{3}, \frac{\xi_{2}}{\xi_{3}})\right
)^{k} B (x', \xi') ,
\end{equation}
where $ B (x', \xi')  $ denotes an analytic $ 2 $-vector defined and
non-zero in a conic neighborhood of $ (0, e_{3}) $. 
The existence of such an integer $ k
$ is a consequence of our analyticity assumption.

Our aim is to draw some consequences from the linearity
of the left hand side of Equation (\ref{5.8}) with respect to $ \xi' $.

Assume first that $ k > 1 $ in (\ref{5.8}). Then taking the $ \xi
$-gradient, we get
$$ 
A(x') = O\left( \left( x_{2} - \chi(x_{3},
    \frac{\xi_{2}}{\xi_{3}})\right)^{k-1} \right),
$$
which implies that $  \chi(x_{3}, \xi_{2}/\xi_{3}) $ actually
depends only on $ x_{3} $. Thus $ B
(x', \xi')  $ is linear with respect to $ \xi' $, so that 
we obtain
\begin{equation}
\label{5.9}
A(x') \xi' = \left( x_{2} - \chi(x_{3})\right)^{k} \tilde{A}(x') \xi',
\end{equation}
where $ \tilde{A}(x') $ denotes another $ 2 \times 2 $ matrix with
real analytic entries.

Let us now assume that $ k =1 $. Equation (\ref{5.8}) becomes 
\begin{equation}
\label{5.10}
A(x') \xi' = \phi (x', \xi')  B (x', \xi') ,
\end{equation}
where $ B $ is a vector-valued symbol of order 0. 
Recall that we are assuming that
$$ 
\phi(0, e_{3}) = 0, \qquad \frac{\partial \phi}{\partial \xi'}(0,
e_{3}) = 0, \qquad {\hbox{ and }}\qquad \frac{\partial \phi}{\partial
x_2}(0, e_{3})
\neq 0.
$$
Since the vanishing in (\ref{5.10}) is of the first order, we have
that $ B(0, e_{3})$  $ \neq 0 $; in particular we may assume that
$$ 
B(0, e_{3}) = (b_{2}(0, e_{3}), b_{3}(0, e_{3}))
$$
and 
\begin{equation}\label{5.11}
b_{3}(0, e_{3})
\neq 0.
\end{equation}
This is no restriction since we can always interchange the second and
the third vector fields.
Taking 
the $ \xi' $-gradient of (\ref{5.10}) and computing everything
at $ (0, e_{3}) $, we easily see that $ A(0) = 0 $. Hence
\begin{equation}
\label{5.12} 
A(x') = x_{2} A^{(2)}(x') + x_{3} A^{(3)}(x'),
\end{equation}
where the $ A^{(j)} $ are real analytic $ 2 \times 2 $ matrices, $ j =
2, 3 $. From this equation we obtain
$$ 
\frac{\partial}{\partial x_{2}} A(x') \xi' = \frac{\partial
  \phi}{\partial x_{2}} B (x', \xi') + \phi (x', \xi') \frac{\partial
  B}{\partial x_{2}} (x', \xi') ,
$$
which, when computed at $ (0, e_{3}) $, yields
$$ 
\frac{\partial A (0)}{\partial x_{2}} 
\begin{bmatrix}
0\\
1
\end{bmatrix}
= \frac{\partial \phi}{\partial x_{2}} (0, e_{3}) B(0, e_{3}).
$$
Let us now consider the second component of the above equation: we
have, from (\ref{5.11}),
\begin{multline*}
\frac{\partial}{\partial x_{2}}
    \left[ \left( x_{2} a^{(2)}_{32}(x') +
    x_{3} a^{(3)}_{32}(x') \right) \xi_{2} \right . \\[5pt]
    + \left . \left . \left( x_{2}
    a^{(2)}_{33}(x') + x_{3} a^{(3)}_{33}(x')\right) \xi_{3} \right ]
\right|_{x'=0\ \xi'=e_{3}} \neq 0,
\end{multline*}
from which we deduce that
\begin{equation}
\label{5.13}
a^{(2)}_{33}(0) \neq 0.
\end{equation}

The second line of the equation $ A(x') \xi' =0 $ then reads:
\begin{equation}
\label{5.14}
\left( x_{2} a^{(2)}_{32}(x') + x_{3} a^{(3)}_{32}(x')\right) \xi_{2}
+ \left( x_{2} a^{(2)}_{33}(x') + x_{3} a^{(3)}_{33}(x')\right)
\xi_{3} = 0.
\end{equation}
Because of (\ref{5.13}) this is the equation of an
analytic submanifold of codimension one containing the point $ (0,
e_{3})
$ and, since $ b_{3}(0, e_{3}) \neq 0$, (\ref{5.10}) implies that
(\ref{5.14}) is equivalent to $ \phi (x', \xi') = 0 $.

Thus we are allowed to change notation and write
\begin{multline}
\label{5.15}
\phi (x', \xi') = \left( x_{2} a^{(2)}_{32}(x') + x_{3}
  a^{(3)}_{32}(x')\right) \xi_{2} \\
+ \left( x_{2} a^{(2)}_{33}(x') + x_{3} a^{(3)}_{33}(x')\right)
\xi_{3},
\end{multline} 
where
\begin{equation}
\label{5.16}
a^{(2)}_{33}(0) \neq 0,
\end{equation}
(recall that we are in Case II${}_{x_2}$ where $\phi_{x_2}\neq 0$ at
$(0,e_3)$).

The following Lemma will help distinguish between two very different
types of families of vector fields. They are both in Case II and will be
denoted Case IIa and Case IIb, which of course will be further
subscripted according to whether $\phi_{x_2}\neq 0$ or
$\phi_{x_3}\neq 0.$
\begin{lemma}
\label{lemma5.1}
Let $ \lambda $ and $ \mu $ be real analytic functions defined in a
neghborhood of the origin and consider the vector field
$$ 
Y = \lambda(x') \frac{\partial}{\partial x_{2}} + \mu(x')
\frac{\partial}{\partial x_{3}}.
$$
Assume that the symbol of $ Y $, $ \lambda(x') \xi_{2} + \mu(x')
\xi_{3} $, vanishes where $ \phi $ vanishes, $ \phi $ being defined in
Equation (\ref{5.15}). Then two cases may occur: 
\begin{itemize}
\item[(a)]{} The set $ \phi^{-1}(0) $ is cylindrical in the $ \xi'
  $-fibers. Then $ \phi (x', \xi') = 0 $ if and only if $ g(x') = 0, $
  for a suitable analytic function $ g $ defined in a neighborhood of
  the origin and having a non-zero $ x' $-gradient. In this case 
$$ 
Y (x', \xi') = g(x') \tilde{Y} (x', \xi') ,
$$
for a suitable vector field $ \tilde{Y} $.
\item[(b)]{} The set $ \phi^{-1}(0) $ is not cylindrical in the $
  \xi'$-fibers. Then there exists an analytic function of $ x' $, $
  h(x') $, defined near $ 0 $, such that
$$ 
Y (x', \xi') = h(x') \phi (x', \xi').
$$
\end{itemize}
\end{lemma}
\begin{proof}
Let us write the function $ \phi $ in (\ref{5.15}) as
$$ 
\phi (x', \xi') = \alpha(x') \xi_{2} + \beta(x') \xi_{3},
$$
where, by (\ref{5.16}), $ \partial \beta (0)/\partial x_{2} \neq 0
$. The vanishing of the symbol of $ Y $ where $ \phi $ vanishes can be
expressed by the following equation:
$$ 
\lambda(x') \xi_{2} + \mu(x') \xi_{3} = a (x', \xi') \left( \alpha(x')
  \xi_{2} + \beta(x') \xi_{3} \right),
$$
where $ a $ is a suitable analytic symbol of order 0---actually
homogenous of degree zero---defined near the
point $ (0, e_{3}) $.

Dividing by $ \xi_{3} $, which is non-zero near $ e_{3} $, and writing
$ \sigma = \xi_{2}/ \xi_{3} $, we have
$$ 
a(x', \sigma, 1) = \frac{\lambda(x') \sigma + \mu(x') }{\alpha(x')
  \sigma + \beta(x')}, \qquad |\sigma| \leq C,
$$
for a suitable positive constant $ C $. Thus
$$ 
a(x', 0, 1) = \frac{\mu(x')}{\beta(x')},
$$
which is also analytic with respect to the variable $ x' $ near the
origin. Since $ \beta(0) = 0 $ and $ \partial \beta (0)/\partial x_{2}
\neq 0$, we have that $ \beta^{-1}(0) $ is a regular analytic curve in
$ \R^{2} $ near the origin. Hence there exists an analytic function $
\gamma(x') $ defined near $ 0 $ such that 
$$ 
\mu(x') = \gamma(x') \beta(x'),
$$
and thus
$$ 
a(x', \sigma, 1) = \frac{\lambda(x') \sigma + \gamma(x')
  \beta(x')}{\alpha(x') \sigma + \beta(x')}.
$$

Now we have

$$ 
\frac{\partial}{\partial \sigma} \frac{\lambda \sigma + \gamma
  \beta}{\alpha \sigma + \beta} = \frac{\beta (\lambda - \alpha
  \gamma)}{(\alpha \sigma + \beta )^{2}}
$$
$$ 
\left( \frac{\partial}{\partial \sigma} \right)^{h+1}
\frac{\lambda\sigma + \gamma \beta }{\alpha \sigma + \beta} =
(-1)^{h} (h+1)! \beta (\lambda - \alpha \gamma)
\frac{\alpha^{h}}{(\alpha \sigma + \beta )^{h+2}}.
$$
Setting $ \sigma = 0 $ in the first line and taking into account the
analyticity of the left hand side, we obtain that
$$ 
\lambda(x') = \alpha(x') \gamma(x') + \delta(x') \beta(x'),
$$
for a suitable analytic function $ \delta $ defined near the
origin. On the other hand, for $ \sigma = 0 $, the second line gives
$$ 
\left. \left( \frac{\partial}{\partial \sigma} \right)^{h+1} a(x',
  \sigma, 1) \right |_{\sigma = 0} = (-1)^{h} (h+1)! \delta(x') \left
  ( \frac{\alpha(x')}{\beta(x')} \right)^{h}.
$$
Now two cases may occur:
\begin{itemize}
\item[i)]{}
Assume that $ \beta $ is a factor of $ \alpha $, i.e. that $
\alpha(x') = \eta(x') \beta(x')^{k} $, for a suitable positive integer
$ k $ and a suitable analytic function $ \eta $. 
In this case $ \phi (x', \xi') = ( \eta \beta^{k-1} \xi_{2} + \xi_{3})
\beta(x')$, with $ \beta(0) = 0 $, $ \partial \beta (0) / \partial
x_{2} \neq 0 $ and $ (\eta \beta^{k-1} \xi_{2} + \xi_{3} )|_{x'=0\
  \xi' = e_{3}} \neq 0 $. We conclude that $ \phi^{-1}(0) =
\beta^{-1}(0) $, or that $ \phi^{-1}(0) $ is the zero set of a
function of $ x' $ only. Moreover in this case we have that $ \lambda
= ( \eta \beta^{k-1} + \delta) \beta $, so that
\begin{multline*} 
\hbox to 2cm{} \lambda(x') \xi_{2} + \mu(x') \xi_{3} = \\
\beta(x') \left[ \left(
    \eta(x') \beta(x')^{k-1} + \delta(x') \right) \xi_{2} + \gamma(x')
  \xi_{3} \right ],
\end{multline*}
which is the conclusion in part (a) of the statement of the Lemma.
\item[ii)]{}
The function $ \beta $ is not a factor of $ \alpha $, i.e. the
quotient $ \alpha / \beta $ is not analytic near $ 0 $. Then
necessarily we must have that $ \delta = 0 $ if $ \beta = 0 $. But
then it is easy to see that there exists a positive integer $ h $ such
that $ \delta / \beta^{h} $ is not analytic near the origin, unless $
\delta \equiv 0 $ in a neighborhood of the origin. Thus
$$ 
\lambda(x') = \alpha(x') \gamma(x'),
$$
at least in a possibly smaller neighborhood of the origin. The above
equation implies that 
$$ 
\lambda(x') \sigma + \mu(x') = \gamma(x') \left( \alpha(x') \sigma +
  \beta(x') \right),
$$
which is the desired conclusion for part (b) of the Lemma.
\end{itemize}
\end{proof}
Summing up we can state the following
\begin{proposition}
\label{proposition5.1}
Assume that the quantity $ A(0, x') \xi' $ vanishes exactly  on an
analytic submanifold $ \Sigma_{p} $ of codimension one inside $
\Sigma_{1} $. Let us denote by $ \phi (x', \xi') = 0 $ a (microlocal)
equation of $ \Sigma_{p} $ near the point $ (0, e_{3}) \in \Sigma_{p}
$. Then the following cases may occur:
\begin{itemize}
\item[I) ]{}
If 
$$ 
\frac{\partial \phi}{\partial \xi'}(0, e_{3}) \neq 0,
$$
then necessarily $  \partial \phi / \partial \xi_{2}(0, e_{3})
\neq 0$ and the equation $ \phi (x', \xi') = 0 $ is equivalent to 
\begin{equation}
\label{5.17:xi_2=0}
\xi_{2} = 0,
\end{equation}
provided a suitable change of coordinates is performed near the $ x'
$-origin. 
\par
In particular we deduce that in this case $ \rank A(0) = 1 $ so that,
on $ \Sigma_{p} $ we also have $ \rank
A(x') =1 $ near the origin.
\item[II)]{}
Assume that 
$$ 
\frac{\partial \phi}{\partial \xi'}(0, e_{3}) = 0,
$$
and 
$$ 
\frac{\partial \phi}{\partial x'}(0, e_{3}) \neq 0.
$$
Then the following cases may occur:
\begin{itemize}
\item[(a)]{}
The equation of $ \Sigma_{p} $ relatively to $ \Sigma_{1} $ does not
depend on $ \xi' $, i.e. $ \Sigma_{p} $ is cylindrical with respect to
the $ \xi' $-fibers. Then we may change coordinates near the origin in
such a way that, in $ \Sigma_{1} $, $ \Sigma_{p} $ is defined by the
equation 
\begin{equation}
\label{5.18:x_j=0}
x_{j} = 0,
\end{equation}
where $ j \in \{ 2, 3\} $.
\item[(b)]{}
Denote by $ \phi (x', \xi') = 0 $ the equation of $ \Sigma_{p} $ in $
\Sigma_{1} $. Then if $ \partial \phi / \partial x_{2} (0, e_{3}) \neq
0$ in a suitable system of coordinates near the origin $ \phi $ is
equivalent to 
\begin{equation}
\label{5.19:b1}
Y(x', \xi') \equiv \lambda(x') \xi_{2} + x_{2} \xi_{3} = 0.
\end{equation}
Here $ \lambda $ denotes a real analytic function such that $
\lambda(0) = 0 $.
\par\noindent
On the other hand assume that $ \partial \phi / \partial x_{2} (0,
e_{3}) = 0 $ and that $ \partial \phi / \partial x_{3} (0, e_{3}) \neq 0
$. Then the equation $ \phi = 0 $ is equivalent to 
\begin{equation}
\label{5.20:b2}
Y(x', \xi') \equiv  \lambda(x') \xi_{2} + \mu(x') \xi_{3} + x_{3}
\xi_{3} = 0,
\end{equation}
where $ \lambda(0) = 0 $, $ \mu(0) = 0 $, $ d_{x'}
\mu (0) = 0
$.
\end{itemize}
\end{itemize}
\end{proposition}
\begin{proof}
To prove the above statement we need only remark that in
Case I any equation of the form $ \xi_{2} - \chi(x') \xi_{3} = 0 $ may
be written as $ \xi_{2} = 0 $, performing a change of coordinates that
leaves $ x_{1} $ unchanged.

As for Case IIa it suffices to notice that $ \Sigma_{p} $ is given, by
what has been shown previously, by the equation $ \beta(x') = 0 $ with
$ d_{x'}\beta (0) \neq 0 $. Thus we can always change coordinates in
the $ (x_{2}, x_{3}) $-plane in such a way that $ \beta(x') = 0 $
becomes $ x_{2} = 0 $ if $ \partial \beta / \partial x_{2} (0) \neq 0
$, or $ x_{3} = 0 $ otherwise.

Let us consider the Case IIb. If $ \partial \phi / \partial x_{2} (0,
e_{3}) \neq 0 $, we have
$$ 
\phi (x', \xi') = \left( x_{2} a^{(2)}_{k2}(x') + x_{3}
  a^{(3)}_{k2}(x') \right)\xi_{2} + \left( x_{2}a^{(2)}_{k3}(x') +
  x_{3} a^{(3)}_{k3}(x') \right) \xi_{3} 
$$
where $ k = 2 $ or $ k = 3 $ depending on which component of the
2-vector $ B $ in (\ref{5.10}) is elliptic at $ (0, e_{3}) $;
moreover $ a^{(2)}_{k3}(0) \neq 0 $. Then we conclude that the equation $
x_{2} a^{(2)}_{k3}(x') + x_{3}a^{(3)}_{k3}(x') = 0 $ is equivalent to
$ x_{2} - \chi(x_{3}) = 0 $, for a suitable analytic function $ \chi $
defined near the origin. Let us perform the following change of
variables in the $ (x_{2}, x_{3}) $-plane:
$$ 
\left\{
\begin{array}{rcl}
y_{2} & = & x_{2} - \chi(x_{3}) \\[7pt]
y_{3} & = & x_{3}
\end{array} \right .
\qquad
\left\{
\begin{array}{rcl}
\eta_{2} & = & \xi_{2} \\
\eta_{3} & = & \xi_{3} + \frac{\textstyle \partial \chi
  (x_{2})}{\textstyle \partial x_{3}} \xi_{2}
\end{array} \right . .
$$
Then in the new coordinates, modulo a non-vanishing factor, we have
$$ 
\phi (x', \xi') = e(x') (\lambda(x') \xi_{2} + x_{2} \xi_{3} ),
$$
which gives (\ref{5.19:b1}).

Assume now that $ \partial \phi / \partial x_{2} (0, e_{3}) = 0 $ and
that $ \partial \phi / \partial x_{3}(0, e_{3}) \neq 0 $. In the above
expression of $ \phi $ we then have $ a^{(3)}_{k3}(0) \neq 0 $ and $
a^{(2)}_{k3}(0) = 0 $, otherwise we would be in the same situation as
above. 

Thus
$$ 
\phi (x', \xi') = a^{(3)}_{k3}(x') \left[ \lambda(x') \xi_{2} +
  \mu(x') \xi_{3} + x_{3} \xi_{3} \right],
$$
with $ \lambda(0) = 0 $ and $ \mu (x') = O(|x'|^{2}) $ which yields
Equation (\ref{5.20:b2}). This completes
the proof of the proposition.
\end{proof}

{\bf Remark.}The seemingly pedantic distinction between the $ x_{2}
$ and $ x_{3} $ variable in the proof above will be useful in
subsequent work, where we shall be concerned with the Gevrey
(analytic) hypoellipticity properties of our operators. The basic tool for
us are microlocal a priori estimates and we shall see that, from a
microlocal point of view, the Gevrey hypoellipticity thresholds for
cases (\ref{5.19:b1}) and (\ref{5.20:b2}), near the same point $ (0,
e_{3}) $, are very different. Naturally, near different
base points, both Cases IIa and IIb may occur for the same operator,
yielding different microlocal hypoellipticity results and the expected
local result.

\vskip 0.4cm
Using Proposition \ref{proposition5.1} we can write the vector fields in a
simpler way.

\begin{proposition}
\label{2-brackets}
The vector fields $ X_{1} $, $ X_{2} $ and $ X_{3} $ satisfying 
hypotheses (A1)-(A4) can be written in the following way:

Case I:
\begin{eqnarray}
\label{5.21:caseI}
X_{1}(x, \xi) & = & \xi_{1} \nonumber\\
X_{2}(x, \xi) & = & a_{21}(x) \xi_{1} + x_{1}^{p-1} \left [\alpha(x')
  \xi_{2} \right. \nonumber \\
& & \left. + x_{1} \left\{
    \tilde{a}_{22}(x)\xi_{2} + \tilde{a}_{23}(x) \xi_{3} \right \}
\right]   \\
X_{3}(x, \xi) & = & a_{31}(x) \xi_{1} + x_{1}^{p-1} \left [ \lambda(x')
  \alpha(x') \xi_{2}
\right. \nonumber \\
&  &  \left . +
x_{1} \left\{
    \tilde{a}_{32}(x)\xi_{2} + \tilde{a}_{33}(x) \xi_{3} \right \}  \right]
   , \nonumber
\end{eqnarray}
for suitable functions $ \alpha(x')  $, with $ \alpha(0) \neq 0 $, and
$ \lambda(x') $. 

Case IIa:
\begin{eqnarray}
\label{5.22:caseIIa}
X_{1}(x, \xi) & = & \xi_{1} \nonumber\\
X_{2}(x, \xi) & = & a_{21}(x) \xi_{1} + x_{1}^{p-1} \left[\, x_{j}
  \left (\tilde{a}_{22}(0, x') 
  \xi_{2} + \tilde{a}_{23}(0, x') \xi_{3} \right) \right. \nonumber \\
& & \left. + x_{1} \left\{
    \hat{a}_{22}(x)\xi_{2} + \hat{a}_{23}(x) \xi_{3} \right \}
\right]   \\
X_{3}(x, \xi) & = & a_{31}(x) \xi_{1} + x_{1}^{p-1} \left [ \,x_{j}
  \left ( \tilde{a}_{32}(0, x') \xi_{2} + \tilde{a}_{33}(0, x') \xi_{3}\right) 
\right. \nonumber \\
&  &  \left . +
x_{1} \left\{
    \hat{a}_{32}(x)\xi_{2} + \hat{a}_{33}(x) \xi_{3} \right \}  \right]
   . \nonumber
\end{eqnarray}
where $ j $ is equal to 2 or 3.

Case IIb:
\begin{eqnarray}
\label{5.23:caseIIb}
X_{1}(x, \xi) & = & \xi_{1} \nonumber\\
X_{2}(x, \xi) & = & a_{21}(x) \xi_{1} + x_{1}^{p-1} \left [ \alpha(x')
  Y (x', \xi')  \right. \nonumber \\
& & \left. + x_{1} \left\{
    \hat{a}_{22}(x)\xi_{2} + \hat{a}_{23}(x) \xi_{3} \right \}
\right]   \\
X_{3}(x, \xi) & = & a_{31}(x) \xi_{1} + x_{1}^{p-1} \left [ 
  \beta(x') Y (x', \xi') 
\right. \nonumber \\
&  &  \left . +
x_{1} \left\{
    \hat{a}_{32}(x)\xi_{2} + \hat{a}_{33}(x) \xi_{3} \right \}  \right]
   . \nonumber
\end{eqnarray}
where $ \alpha $ is a non-vanishing analytic function defined in a
neighborhood of the origin, $ \beta $ is analytic and $ Y (x', \xi')
$ is a vector field of the form (\ref{5.19:b1}) or (\ref{5.20:b2}).
\end{proposition}
\begin{proof}
Case I is straightforward, due to Proposition \ref{proposition5.1}. The same
proposition also implies Case IIa. Case IIb follows from Proposition
\ref{proposition5.1} and Lemma \ref{lemma5.1} (b).
\end{proof}
{\bf Remark.}
We point out that, since $ \Sigma_{1} = \{ x_{1} = 0, \;\xi_{1} = 0 \}
$, the forms (\ref{5.21:caseI}) - (\ref{5.23:caseIIb}) for our vector fields
actually have some further properties, which will turn out to be
important for the regularity estimates. Basically these properties
state that the fields are linearly independent outside of the
characteristic manifold and that the number of layers of the Poisson
stratification is finite. We postpone a precise statement of this fact
until the final step in order not to burden the exposition too much.

\vskip 0.4cm

The next step consists in using Assumption (A3) and the remaining part
of (A4) to make the form of the vector fields more precise.
\section{Finer forms for the vector fields}
\setcounter{equation}{0}
\setcounter{theorem}{0}
\setcounter{proposition}{0}  
\setcounter{lemma}{0}
\setcounter{corollary}{0} 
\setcounter{definition}{0}

\subsection{Case I}
\setcounter{equation}{0}
\setcounter{theorem}{0}
\setcounter{proposition}{0}  
\setcounter{lemma}{0}
\setcounter{corollary}{0} 
\setcounter{definition}{0}
\renewcommand{\thetheorem}{\thesubsection.\arabic{theorem}}
\renewcommand{\theproposition}{\thesubsection.\arabic{proposition}}
\renewcommand{\thelemma}{\thesubsection.\arabic{lemma}}
\renewcommand{\thedefinition}{\thesubsection.\arabic{definition}}
\renewcommand{\thecorollary}{\thesubsection.\arabic{corollary}}
\renewcommand{\theequation}{\thesubsection.\arabic{equation}}

By Proposition \ref{2-brackets} we are dealing with the fields:
\begin{eqnarray*}
X_{1}(x, \xi) & = & \xi_{1} \\[3pt]
\begin{bmatrix}
X_{2}(x, \xi) \\ X_{3}(x, \xi) 
\end{bmatrix}
& = & 
\begin{bmatrix}
  a_{21}(x) \\ a_{31}(x)
\end{bmatrix}
\xi_{1} + x_{1}^{p-1} \left\{
\begin{bmatrix}
  \alpha(x') & 0 \\
  \lambda(x') \alpha(x') & 0
\end{bmatrix}
\xi' + x_{1} \tilde{A}(x) \xi' \right\},
\end{eqnarray*}
with obvious notation. We can see at once that the only brackets
that matter are
$$ 
\ad^{j}(X_{1}) X_{k}, \qquad k = 2, 3, \qquad j= p, p+1, \ldots ,q-2.
$$
The above quantity vanishes on $ \Sigma_{p} = \{ x_{1} = \xi_{1}
= 0, \xi_{2} = 0 \}$, so that, taking $ j = p $, we conclude that
$$ 
\tilde{A}(x') 
\begin{bmatrix}
0 \\
\xi_{3}
\end{bmatrix}
= 0,
$$
which implies that
$$ 
\tilde{a}_{2 3}(x){\Big |_{x_{1} = 0}} = \tilde{a}_{33}(x){\Big
  |_{x_{1} = 0}} = 0.
$$
Thus we may write the fields as
\begin{eqnarray*}
X_{1}(x, \xi) & = & \xi_{1} \nonumber \\[3pt]
\begin{bmatrix}
  X_{2}(x, \xi)  \\ X_{3}(x, \xi) 
\end{bmatrix}
& = &
\begin{bmatrix}
  a_{21}(x) \\ a_{31}(x) 
\end{bmatrix}
\xi_{1} + x_{1}^{p-1} \left\{
\begin{bmatrix}
  \alpha(x') + x_{1} \tilde{a}_{22}(x)  & 0 \\
  \lambda(x') \alpha(x') + x_{1}\tilde{a}_{32}(x) & 0
\end{bmatrix}
\xi' \right . \nonumber \\ 
& & \hskip 3cm \left . + x_{1} 
\begin{bmatrix}
  \tilde{a}_{23}(x) \\ \tilde{a}_{33}(x)
\end{bmatrix}
\xi_{3} \right\},
\end{eqnarray*}
for suitable analytic coefficients $ \tilde{a}_{i 3} $, $ i = 1, 2 $. 

Proceeding analogously and using the remaining brackets, we
conclude that
\begin{eqnarray}
\label{6.1.1}
X_{1}(x, \xi) & = & \xi_{1} \nonumber \\[3pt]
\begin{bmatrix}
  X_{2}(x, \xi)  \\ X_{3}(x, \xi) 
\end{bmatrix}
& = &
\begin{bmatrix}
  a_{21}(x) \\ a_{31}(x) 
\end{bmatrix}
\xi_{1} + x_{1}^{p-1} \left\{
\begin{bmatrix}
  \alpha(x') + x_{1} \tilde{a}_{22}(x)  & 0 \\
  \lambda(x') \alpha(x') + x_{1}\tilde{a}_{32}(x) & 0
\end{bmatrix}
\xi' \right . \nonumber \\ 
& & \hskip 3cm \left . + x_{1}^{q-p} 
\begin{bmatrix}
  \tilde{a}_{23}(x) \\ \tilde{a}_{33}(x)
\end{bmatrix}
\xi_{3} \right\},
\end{eqnarray}
for suitable analytic functions $ \alpha(x') \neq 0 $ (as always in Case
I - cf.(\ref{5.21:caseI})), $\tilde{a}_{i 3}$,
$ i = 1, 2 $, and $
\lambda(x'). $

Furthermore the ellipticity of the Poisson brackets of length $ q $
tells us that 
\begin{equation}
\label{6.1.2}
\tilde{A}(0, x') 
\begin{bmatrix}
0\\
1
\end{bmatrix}
\neq 0.
\end{equation}
On the other hand, the fields $ X_{2} $, $ X_{3}
$ in (\ref{6.1.1}) are linearly independent for $ x_{1} \neq 0 $ if and
only if
$$ 
x_{1}^{q-p} \det
\begin{bmatrix}
\alpha(x') + x_{1} \tilde{a}_{22}(x) & \tilde{a}_{23}(x) \\[7pt]
\lambda(x') \alpha(x') + x_{1} \tilde{a}_{23}(x) & \tilde{a}_{33}(x)
\end{bmatrix}
\neq 0,
$$
i.e.
$$ 
- \lambda(x') \tilde{a}_{23}(x) + \tilde{a}_{33}(x) +
\frac{x_{1}}{\alpha(x')} \det \tilde{A}(x) \neq 0,
$$
if $ x_{1} \neq 0 $, or
\begin{equation}
\label{6.1.3}
\det \left(
\begin{bmatrix}
1 \\
\lambda(x') 
\end{bmatrix} \otimes
\begin{bmatrix}
1\\
0
\end{bmatrix} + 
\frac{x_{1}}{\alpha(x')} \tilde{A}(x) \right) \neq 0,
\end{equation}
if $ x_{1} \neq 0 $. Another way of stating the above condition is 
\begin{equation}
\label{6.1.4}
\langle \begin{bmatrix}
- \lambda (x') \\
1
\end{bmatrix} , \tilde{A}(x) 
\begin{bmatrix}
0\\
1
\end{bmatrix} \rangle +
\frac{x_{1}}{\alpha(x')} \det \tilde{A}(x) \neq 0,
\end{equation}
if $ x_{1} \neq 0 $.

\subsection{Case IIa}
\setcounter{equation}{0}
\setcounter{theorem}{0}
\setcounter{proposition}{0}  
\setcounter{lemma}{0}
\setcounter{corollary}{0} 
\setcounter{definition}{0}
\renewcommand{\thetheorem}{\thesubsection.\arabic{theorem}}
\renewcommand{\theproposition}{\thesubsection.\arabic{proposition}}
\renewcommand{\thelemma}{\thesubsection.\arabic{lemma}}
\renewcommand{\thedefinition}{\thesubsection.\arabic{definition}}
\renewcommand{\thecorollary}{\thesubsection.\arabic{corollary}}
\renewcommand{\theequation}{\thesubsection.\arabic{equation}}

We begin by considering the fields in (\ref{5.22:caseIIa}) and
again use Assumption (A4) and (A3). Thanks to the remarks made
above, we can see that, taking $ p $ derivatives with respect to $ x_{1}
$, we have
$$ 
\hat{A} (x) \xi' = 0 \quad \text{if }\quad x_{1}
= x_{j} = 0, \quad j= 2,3,
$$
i.e.
\begin{equation}
\label{6.2.1}
\hat{A}(x) = x_{1} \hat{A}_{1}(x) + x_{j} \hat{A}_{2}(x).
\end{equation}
Hence $ X_{2} $ and $ X_{3} $ can be written:
\begin{multline*} 
\begin{bmatrix}
  X_{2}(x, \xi)  \\ X_{3}(x, \xi) 
\end{bmatrix}
=
\begin{bmatrix}
  a_{21}(x) \\ a_{31}(x)
\end{bmatrix}
\xi_{1} + x_{1}^{p-1} \left\{ x_{j} \tilde{A}(x') \xi' \right .\\
\left . + \left( x_{1}^{2} \hat{A}_{1}(x) + x_{1} x_{j}
  \hat{A}_{j}(x)\right)\xi' \right \}
\end{multline*}
$$ 
= \begin{bmatrix}
  a_{21}(x)  \\ a_{31}(x) 
\end{bmatrix}
\xi_{1} + x_{1}^{p-1} \left\{ x_{j} \tilde{\tilde{A}}(x) \xi' +
  x_{1}^{2} \hat{A}(x)\xi'\right\},
$$
the meaning of the symbols being obvious.

Iterating this argument we reach the following form for the vector
fields: 
\begin{eqnarray}
\label{6.2.2}
X_{1}(x, \xi) & = & \xi_{1} \nonumber \\
X_{2}(x, \xi) & = & a_{21}(x) \xi_{1} + x_{1}^{p-1} \left[ x_{j} \left(
    \tilde{a}_{22}(x) \xi_{2} + \tilde{a}_{23}(x) \xi_{3} \right)
\right . \nonumber \\
& & \left . + x_{1}^{q-p} \left(\hat{a}_{22}(x) \xi_{2} + \hat{a}_{23}(x)
    \xi_{3} \right) \right] \\
X_{3}(x, \xi) & = & a_{31}(x) \xi_{1} + x_{1}^{p-1} \left[ x_{j} \left(
\tilde{a}_{32}(x) \xi_{2} + \tilde{a}_{33}(x) \xi_{3} \right)
\right. \nonumber \\
& & \left. + x_{1}^{q-p} \left( \hat{a}_{32}(x) \xi_{2} + \hat{a}_{33}(x)
    \xi_{3} \right) \right], \nonumber 
\end{eqnarray}
$ j \in \{2, 3\}. $

Proceeding as in Case I we see that the
ellipticity of the last Poisson layer means that
\begin{equation}
\label{6.2.3:ellcase2}
\det \hat{A}(x) \Big|_{x_{1}=0, \ x_{j} = 0} \neq 0.
\end{equation}
On the other hand, Assumption (A3) together with
(\ref{6.2.3:ellcase2}) means that
\begin{equation}
\label{6.2.4:AssA3case2}
\det \left( x_{j}\tilde{A}(x) + x_{1}^{q-p} \hat{A}(x)\right) \neq 0,
\end{equation}
if $ x_{1} \neq 0 $.

\subsection{Case IIb}
\setcounter{equation}{0}
\setcounter{theorem}{0}
\setcounter{proposition}{0}  
\setcounter{lemma}{0}
\setcounter{corollary}{0} 
\setcounter{definition}{0}
\renewcommand{\thetheorem}{\thesubsection.\arabic{theorem}}
\renewcommand{\theproposition}{\thesubsection.\arabic{proposition}}
\renewcommand{\thelemma}{\thesubsection.\arabic{lemma}}
\renewcommand{\thedefinition}{\thesubsection.\arabic{definition}}
\renewcommand{\thecorollary}{\thesubsection.\arabic{corollary}}
\renewcommand{\theequation}{\thesubsection.\arabic{equation}}

Let us consider the fields in (\ref{5.23:caseIIb}) and use Assumptions (A3)
and (A4). Taking the $ p $-th derivative with respect to $ x_{1} $ we
obtain that 
$$ 
\hat{A}(x) \xi' = 0 \qquad \text{if} \quad x_{1}=0 \quad \text{and}
\quad Y (x', \xi') = 0.
$$
By Lemma \ref{lemma5.1} (b), this implies that there is an analytic
2-vector, $h^{(1)}(x'),$ defined near the origin, such that
$$ 
\hat{A}(0, x') \xi' = h^{(1)}(x') Y (x', \xi') , \qquad h^{(1)}(x') = 
\begin{bmatrix}
  h^{(1)}_{2}(x') \\ h^{(1)}_{3}(x')
\end{bmatrix},
$$
so that
$$ 
\hat{A}(x) \xi' = h^{(1)}(x') Y (x', \xi')  + x_{1} \hat{A}^{(1)}(x) \xi'
$$
and hence
\begin{multline*}
\begin{bmatrix}
  X_{2}(x, \xi)  \\ X_{3}(x, \xi) 
\end{bmatrix}
=
\begin{bmatrix}
  a_{21}(x) \\ a_{31}(x)
\end{bmatrix}
\xi_{1} + \\x_{1}^{p-1} \left\{ \left( 
\begin{bmatrix}
  \alpha(x') \\  \beta(x')
\end{bmatrix}
+ x_{1} h^{(1)}(x') \right) Y (x', \xi') \right .\\
  + x_{1}^{2} \hat{A}^{(1)}(x) \xi' \Big \}
\end{multline*}
Iterating this argument we obtain that
\begin{eqnarray}
\label{6.3.1}
X_{1}(x, \xi) & = & \xi_{1} \nonumber \\[3pt]
\begin{bmatrix}
  X_{2}(x, \xi)  \\ X_{3}(x, \xi) 
\end{bmatrix}
& = &
\begin{bmatrix}
  a_{21}(x) \\ a_{31}(x)
\end{bmatrix}
\xi_{1} \\
&  &
+ x_{1}^{p-1} \left\{ h(x) Y (x', \xi') + x_{1}^{q-p}
  \hat{A}(x) \xi' \right\}, \nonumber
\end{eqnarray}
where $ h(x) $ is a 2-vector function, $ h(x) = ( h_{2}(x), h_{3}(x))
$, such that $ h_{2}(0) \neq 0 $, and $ \hat{A} $ is a $ 2 \times 2 $
matrix with real analytic entries defined near the origin.

Assumption (A4) then implies that $ \hat{A}(x)\xi' $ cannot be zero if
$ x_{1} = 0 $, $ \xi_{1} = 0 $ and $ Y (x', \xi') = 0 $; but, since $
Y(0, \xi') \equiv 0 $ for every $ \xi' \in \R^{2} $, 
we easily get that
\begin{equation}
\label{6.3.2}
\det \hat{A}(0) \neq 0,
\end{equation}
while the linear independence of the vector fields outside of $
\Sigma_{1} $ yields
\begin{equation}
\label{6.3.3:linindc2b}
\det ( h(x) \otimes Y (x') + x_{1}^{q-p} \hat{A}(x)) \neq 0,
\end{equation}
if $ x_{1} \neq 0 $. Here $ Y(x') $ denotes the 2-vector whose
components are the coefficients of the vector field $ Y $.

We summarize the above argument in
\begin{theorem}
\label{th1}
Let $ X_{1}$, $X_{2}$, $X_{3} $ satisfy Assumptions (A1) - (A4). Then
there is a suitable system of coordinates defined in a neighborhood of
the point $ (0, e_{3}) $, such that the field can be written in one of
the following ways:
\begin{itemize}
\item[Case I)]{}
\begin{eqnarray}
\label{6.3.4}
X_{1}(x, \xi) & = & \xi_{1} \nonumber \\
X_{2}(x, \xi) & = & a_{21}(x) \xi_{1} + x_{1}^{p-1} \left[ \left(
    \alpha(x') + x_{1} \tilde{a}_{22}(x) \right)\xi_{2} \right .\\
& & \left . + x_{1}^{q-p}
  \tilde{a}_{23}(x) \xi_{3} \right] \nonumber \\
X_{3}(x, \xi) & = & a_{31}(x) \xi_{1} + x_{1}^{p-1} \left[ \left(
\lambda(x') \alpha(x') + x_{1} \tilde{a}_{32}(x)\right) \xi_{2} \right
. \nonumber \\
& & \left . + x_{1}^{q-p} \tilde{a}_{33}(x) \xi_{3} \right],
\nonumber
\end{eqnarray}
for suitable analytic functions $ \tilde{a}_{ij} $, $ i, j = 2, 3 $, $
\lambda(x') $, and $ \alpha(x') \neq 0 $. Moreover we have
\begin{equation}
\label{6.3.5:caseIlastlayer}
\begin{bmatrix}
  \tilde{a}_{23}(0, x') \\ \tilde{a}_{33}(0, x')
\end{bmatrix}
\neq 0,
\end{equation}
and
\begin{equation}
\label{6.3.6:case1linind}
- \lambda(x') \tilde{a}_{23}(x) + \tilde{a}_{33}(x) +
  \frac{x_{1}}{\alpha(x')} \det \tilde{A}(x) \neq 0,
\end{equation}
if $ x_{1} \neq 0 $.
\item[Case IIa)]{}
\begin{eqnarray}
\label{6.3.7}
X_{1}(x, \xi)  & = & \xi_{1} \nonumber \\
X_{2}(x, \xi)  & = & a_{21}(x)\xi_{1} + x_{1}^{p-1} \left[ x_{j}
  \left( \tilde{a}_{22}(x)\xi_{2} + \tilde{a}_{23}(x) \xi_{3} \right)
\right . \nonumber \\
&  &  \left . + x_{1}^{q-p} \left(\hat{a}_{22}(x)\xi_{2} +
    \hat{a}_{23}(x)\xi_{3} \right) \right] \\
X_{3}(x, \xi) & = & a_{31}(x)\xi_{1} + x_{1}^{p-1} \left[ x_{j} \left
    ( \tilde{a}_{32}(x)\xi_{2} + \tilde{a}_{33}(x)\xi_{3} \right)
\right . \nonumber \\
& & \left . + x_{1}^{q-p} \left( \hat{a}_{32}(x)\xi_{2} +
    \hat{a}_{33}(x)\xi_{3} \right) \right ], \nonumber 
\end{eqnarray}
where $ j \in \{2, 3 \} $, $ \tilde{a}_{ij} $, $ \hat{a}_{ij} $ are
analytic functions, $ i, j = 2, 3 $, such that
\begin{equation}
\label{6.3.8:case2alastlayer}
\det \tilde{A}(x)\Big |_{x_{1}=0} \neq 0, \qquad 
\det \hat{A}(x) \Big |_{x_{1} = x_{j} = 0} \neq 0.
\end{equation}
Moreover
\begin{equation}
\label{6.3.9:case2alinind}
\det \left( x_{j} \tilde{A}(x) + x_{1}^{q-p} \hat{A} (x) \right) \neq 0,
\end{equation}
if $ x_{1} \neq 0 $.
\item[Case IIb)]{}
\begin{eqnarray}
\label{6.3.10}
X_{1}(x, \xi) & = & \xi_{1} \nonumber \\
X_{2}(x, \xi) & = & a_{21}(x) \xi_{1} + x_{1}^{p-1} \left[ h_{2}(x)
  \left ( \alpha(x') \xi_{2} + \beta(x')\xi_{3} \right) \right
. \nonumber \\
& &
\left . + x_{1}^{q-p} \left( \hat{a}_{22}(x)\xi_{2} + \hat{a}_{23}(x)
    \xi_{3} \right) \right] \\
X_{3}(x, \xi) & = & a_{31}(x)\xi_{1} + x_{1}^{p-1} \left[ h_{3}(x)
  \left (\alpha(x') \xi_{2} + \beta(x') \xi_{3}\right) \right
. \nonumber \\
& & 
\left . + x_{1}^{q-p} \left( \hat{a}_{32}(x)\xi_{2} +
    \hat{a}_{33}(x)\xi_{3} \right) \right] \nonumber
\end{eqnarray}
where we may assume that $ h_{2}(0) \neq 0 $, $ h_{j} $ and $
\hat{a}_{ij} $ are suitable analytic functions, and the field $
\alpha(x')\xi_{2} + \beta(x')\xi_{3} $ has the form in (\ref{5.19:b1})
or (\ref{5.20:b2}). Moreover
\begin{equation}
\label{6.3.11:case2blastlayer}
\det \hat{A}(0) \neq 0
\end{equation}
and
\begin{equation}
\label{6.3.12:case2blinind}
\det \left( h(x) \otimes \begin{bmatrix}
  \alpha(x') \\
   \beta(x')
\end{bmatrix}
+ x_{1}^{q-p}\hat{A}(x) \right) \neq 0,
\end{equation}
if $ x_{1} \neq 0 $.
\end{itemize}
\end{theorem}

\section{Examples}
\setcounter{equation}{0}
\setcounter{theorem}{0}
\setcounter{proposition}{0}  
\setcounter{lemma}{0}
\setcounter{corollary}{0} 
\setcounter{definition}{0}
\renewcommand{\thetheorem}{\thesection.\arabic{theorem}}
\renewcommand{\theproposition}{\thesection.\arabic{proposition}}
\renewcommand{\thelemma}{\thesection.\arabic{lemma}}
\renewcommand{\thedefinition}{\thesection.\arabic{definition}}
\renewcommand{\thecorollary}{\thesection.\arabic{corollary}}
\renewcommand{\theequation}{\thesection.\arabic{equation}}

We collect in this section a few examples of the fields obtained in
Theorem \ref{th1}. The Case I examples all have the following
stratification: 
\begin{eqnarray*}
\Sigma_{1} &=& \{ x_{1} = \xi_{1} = 0 \} \\
\Sigma_{2} &=& \Sigma_{1} \\
           & \vdots &  \\
\Sigma_{p} &=& \{ x_{1} = \xi_{1} = 0,\ \xi_{2}=0 \} \\
\Sigma_{p+1}&=& \Sigma_{p} \\
           & \vdots & \\
\Sigma_{q} &=& \{0\},
\end{eqnarray*}
where $ \{0\} $ denotes the zero section of the cotangent bundle.
\subsection{Case I}
\setcounter{equation}{0}
\setcounter{theorem}{0}
\setcounter{proposition}{0}  
\setcounter{lemma}{0}
\setcounter{corollary}{0} 
\setcounter{definition}{0}
\renewcommand{\thetheorem}{\thesubsection.\arabic{theorem}}
\renewcommand{\theproposition}{\thesubsection.\arabic{proposition}}
\renewcommand{\thelemma}{\thesubsection.\arabic{lemma}}
\renewcommand{\thedefinition}{\thesubsection.\arabic{definition}}
\renewcommand{\thecorollary}{\thesubsection.\arabic{corollary}}
\renewcommand{\theequation}{\thesubsection.\arabic{equation}}
\begin{itemize}
\item
Let $ \alpha \equiv 1 $, $ \lambda = 0 $, $ a_{21} = a_{31} = 0 $ and 
$$ 
\tilde{A} = 
\begin{bmatrix}
  0 & 1 \\
  0 & 1
\end{bmatrix}
.
$$
Then we have the fields
$$ 
\xi_{1}, \qquad x_{1}^{p-1} \left[ \xi_{2} + x_{1}^{q-p}
  \xi_{3}\right], \qquad x_{1}^{q-1}\xi_{3}.
$$
\item
Let $ \alpha \equiv 1 $, $ \lambda = 0 $, $ a_{21} = a_{31} = 0 $ and 
$$ 
\tilde{A} = 
\begin{bmatrix}
  0 & 0 \\
  0 & 1
\end{bmatrix}
.
$$
Then we have the fields
$$ 
\xi_{1}, \qquad x_{1}^{p-1}\xi_{2}, \qquad x_{1}^{q-1}\xi_{3},
$$
which is the Oleinik-Radkevi\v c operator.

\item
Let $ \alpha \equiv 1 $, $ \lambda = 0 $, $ a_{21} = a_{31} = 0 $ and 
$$ 
\tilde{A} = 
\begin{bmatrix}
  0 & 1 \\
  x_{1}^{q-p-1} & 0
\end{bmatrix}
.
$$
Then we have the fields
$$ 
\xi_{1}, \qquad x_{1}^{p-1}\left[\xi_{2} + x_{1}^{q-p}\xi_{3}\right],
\qquad x_{1}^{q-1}\xi_{2}.
$$
\end{itemize}
Concerning the conditions of Theorem \ref{th1} we see that
the vector $ (\tilde{a}_{23}, \tilde{a}_{33}) $ is equal to $ (1, 1) $ in the
first case, $ (0, 1) $ in the second case and to $ (1, 0) $ in the third
case. Moreover (\ref{6.3.6:case1linind}) reads as $ 1 + x_{1} \cdot 0
\neq 0 $ in the first and second cases, $ x_{1} \det \tilde{A} = -
x_{1}^{q-p} \neq 0 $ if $ x_{1} \neq 0 $ in the third case. 

\subsection{Case IIa} 
\setcounter{equation}{0}
\setcounter{theorem}{0}
\setcounter{proposition}{0}  
\setcounter{lemma}{0}
\setcounter{corollary}{0} 
\setcounter{definition}{0}
\renewcommand{\thetheorem}{\thesubsection.\arabic{theorem}}
\renewcommand{\theproposition}{\thesubsection.\arabic{proposition}}
\renewcommand{\thelemma}{\thesubsection.\arabic{lemma}}
\renewcommand{\thedefinition}{\thesubsection.\arabic{definition}}
\renewcommand{\thecorollary}{\thesubsection.\arabic{corollary}}
\renewcommand{\theequation}{\thesubsection.\arabic{equation}}
For the Case IIa, the stratification is as for Case I except that
$\Sigma_p$ is now given by: 
$$\Sigma_{p} = \{ x_{1} = \xi_{1} = 0,\  x_{2}=0 \}.
$$
Let us take
$ j = 2
$,
$ a_{21} = a_{31} = 0
$ and
$
\tilde{A} = Id
$. Then from the condition $ \det \left( x_{2} Id + x_{1}^{q-p}
\hat{A} \right) \neq 0
$ if $ x_{1} \neq 0 $ we easily deduce that the matrix $ \hat{A} $ must have
non-zero strictly complex eigenvalues. Set
$$ 
\hat{A} =
\begin{bmatrix}
  \lambda & \mu \\
  - \mu & \lambda
\end{bmatrix}
, \qquad \mu \neq 0.
$$
Then our conditions are satisfied and we obtain the fields
\begin{multline*}
\xi_{1}, \qquad x_{1}^{p-1}\left[ x_{2} \xi_{2} + x_{1}^{q-p} (\lambda \xi_{2}
  + \mu \xi_{3}) \right], \\
x_{1}^{p-1}\left[ x_{2} \xi_{3} + x_{1}^{q-p} (-\mu \xi_{2}
  + \lambda \xi_{3}) \right] .
\end{multline*}

\subsection{Case IIb}
\setcounter{equation}{0}
\setcounter{theorem}{0}
\setcounter{proposition}{0}  
\setcounter{lemma}{0}
\setcounter{corollary}{0} 
\setcounter{definition}{0}
\renewcommand{\thetheorem}{\thesubsection.\arabic{theorem}}
\renewcommand{\theproposition}{\thesubsection.\arabic{proposition}}
\renewcommand{\thelemma}{\thesubsection.\arabic{lemma}}
\renewcommand{\thedefinition}{\thesubsection.\arabic{definition}}
\renewcommand{\thecorollary}{\thesubsection.\arabic{corollary}}
\renewcommand{\theequation}{\thesubsection.\arabic{equation}}
Here the non-symplectic layer $\Sigma_p$ is given near
$(0,\xi_3)$ by: 
$$\Sigma_{p} = x_3^2\xi_2+x_2\xi_3 = 0.
$$
Let $ \phi (x',
\xi') =
\alpha(x')
\xi_{2} +
\beta(x')
\xi_{3} =
\lambda(x') \xi_{2} + x_{2}\xi_{3} $, with $ \lambda 
\not\equiv 0 $, $ \lambda(0) = 0 $, as e.g. in (\ref{5.19:b1}); we may
assume that $ \lambda(x') / x_{2} $ 
is not an analytic function near the origin.

Moreover let $ a_{21} = a_{31} = 0 $, $ h_{2} = 1 $ and $ h_{3} = 0 $. Then we
have the fields 
\begin{multline*}
\xi_{1}, \quad x_{1}^{p-1} \left[ \lambda \xi_{2} + x_{2} \xi_{3} +
  x_{1}^{q-p} \left(\hat{a}_{22}\xi_{2} + \hat{a}_{23} \xi_{3}\right) \right], \quad
x_{1}^{q-1} \left[ \hat{a}_{32} \xi_{2} + \hat{a}_{33}\xi_{3} \right].
\end{multline*}
Conditions (\ref{6.3.11:case2blastlayer}) and (\ref{6.3.12:case2blinind})
become $ \det  \hat{A} \neq 0 $ and $ 
\lambda \hat{a}_{33} - \hat{a}_{32} x_{2} + x_{1}^{q-p} \det \hat{A}$ 
$ \neq 0 $
if $ x_{1} \neq 0 $. If $ q-p $ is e.g. even we may choose $ \lambda =
x_{3}^{2} $, $ \hat{a}_{33} = \sgn \det \hat{A} $, $ \hat{a}_{32} = 0 $ to
write a particular case of the above fields.

\section{The behavior of the bicharacteristic curves and a finer
classification} 
\setcounter{equation}{0}
\setcounter{theorem}{0}
\setcounter{proposition}{0}  
\setcounter{lemma}{0}
\setcounter{corollary}{0} 
\setcounter{definition}{0}
\renewcommand{\thetheorem}{\thesection.\arabic{theorem}}
\renewcommand{\theproposition}{\thesection.\arabic{proposition}}
\renewcommand{\thelemma}{\thesection.\arabic{lemma}}
\renewcommand{\thedefinition}{\thesection.\arabic{definition}}
\renewcommand{\thecorollary}{\thesection.\arabic{corollary}}
\renewcommand{\theequation}{\thesection.\arabic{equation}}
%

%

In this Section we present a classification of the various instances
of the ``sums of squares operators'' in which we get in Case I.

Consider (\ref{6.3.4}); $ X_{1} $ actually denotes the only non
characteristic vector field. Let us consider the null bicharacteristic
curves of $ X_{1} $, $ \gamma_{(\bar{x}, \bar{\xi})}(t) = (\bar{x},
\bar{\xi}) + t (e_{1}, 0) $, where $ \bar{\xi}_{1} = 0 $. If $
\bar{x}_{1} = 0 $, then $ \gamma_{(\bar{x}, \bar{\xi})} (t)  =
\gamma_{(\bar{x}', \bar{\xi}')} (t) = (0, \bar{x}', 0, \bar{\xi}') + t
(e_{1}, 0) $ is actually a null bicharacteristic curve of $ X_{1} $
issued from a point $ (0, \bar{x}', 0, \bar{\xi}') $ of $ \Sigma_{1}
$. Assume $ t \neq 0 $ and compute $ X_{2} $ and $ X_{3} $ on such a
curve. We obtain
\begin{eqnarray}
\label{8.1:x2x3}
X_{2}(\gamma_{(\bar{x}', \bar{\xi}')} (t)) & = & 
t^{p-1} \left[ \left( \alpha(\bar{x}') + t \tilde{a}_{22}(t, \bar{x}')
  \right) \bar{\xi}_{2} \right . \nonumber \\
& & \left .
+ t^{q-p} \tilde{a}_{23}(t, \bar{x}') \bar{\xi}_{3} \right] \\
X_{3}(\gamma_{(\bar{x}', \bar{\xi}')} (t)) & = &
t^{p-1} \left[ \left( \lambda(\bar{x}') \alpha(\bar{x}') + t
    \tilde{a}_{32}(t, \bar{x}') \right) \bar{\xi}_{2} \right
  . \nonumber \\
&  & \left .
+ t^{q-p} \tilde{a}_{33}(t, \bar{x}') \bar{\xi}_{3} \right]. \nonumber
\end{eqnarray}
Assume that the point $ (0, \bar{x}', 0, \bar{\xi}') $ is in a
neighborhood of $ (0, e_{3}) $. Then $ \xi_{3} \neq 0 $ and also $
\alpha(\bar{x}') \neq 0 $ by Theorem \ref{th1}. 
On the other hand nothing is known a priori about the function $
\lambda $. We point out explicitly that we chose $ X_{2} $ as the
field having a non-zero $ \partial / \partial x_{2} $ coefficient near
$ (0, e_{3}) $, thus breaking the $ X_{2} $--$ X_{3} $ symmetry. This
is evidently no restriction of generality, provided we bear in mind
that analogous statements hold if we interchange the roles of $ X_{2}
$ and $ X_{3} $.

When $ t \neq 0 $ we may consider  the characteristic set of $
X_{2}(\gamma_{(\bar{x}', \bar{\xi}')} (t)) $; we obtain that $
X_{2}(\gamma_{(\bar{x}', \bar{\xi}')} (t) = 0 $ if and only if
$$ 
\bar{\xi}_{2} = - t^{q-p} \frac{\tilde{a}_{23}(t,
  \bar{x}')}{\alpha(\bar{x}') + t \tilde{a}_{22}(t, \bar{x}')}
\bar{\xi}_{3}. 
$$
Let us now compute $ X_{3}(\gamma_{(\bar{x}', \bar{\xi}')}(t) \Big
|_{X_{2}(\gamma_{(\bar{x}', \bar{\xi}')}(t)) = 0} $; we get
\begin{multline}
\label{8.2:x3onx20}
\frac{\alpha(\bar{x}')}{\alpha(\bar{x}') + t \tilde{a}_{22}(t,
  \bar{x}')} t^{q-1} \Big[ - \lambda(\bar{x}') \tilde{a}_{23}(t,
  \bar{x}') + \tilde{a}_{33}(t, \bar{x}')  \\
+ \frac{t}{\alpha(\bar{x}')}
  \det \tilde{A}(t, \bar{x}') \Big] \bar{\xi}_{3}
\end{multline}
where the quantity in square brackets is that playing a role in
Equation (\ref{6.3.6:case1linind}) and is non-zero provided $ t \neq 0 $. We
also point out that the coefficient $ \alpha(\bar{x}')
(\alpha(\bar{x}') + t \tilde{a}_{22}(t, \bar{x}') )^{-1} $ is also non-zero at $ t = 0 $.

The above discussion motivates the following

\begin{definition}
\label{8.4:def1r}
We say that the fields $ X_{1} $, $ X_{2} $, $ X_{3} $ of (\ref{6.3.4})
are in Case $ I_{0} $ or of type $ I_{0} $ if
\begin{equation}
\label{8.3:1zero}
- \lambda(0) \tilde{a}_{23}(0, 0) + \tilde{a}_{33}(0, 0) \neq 0.
\end{equation}
This means that, as $ t \to 0 $
$$ 
X_{3}(\gamma_{(\bar{x}', \bar{\xi}')}(t))
\Big|_{X_{2}(\gamma_{(\bar{x}', \bar{\xi}')} (t)) = 0} \sim t^{q-1}, 
$$
uniformly with respect to $ \bar{x}' $, $ \bar{\xi}_{3} \neq 0 $.

Assume now that (\ref{8.3:1zero}) no longer holds and let 
\begin{equation}
\label{1r}
-\lambda(0) \tilde{a}_{23}(t, 0) + \tilde{a}_{33}(t, 0) +
 \frac{t}{\alpha(0)} \det \tilde{A}(t, 0) \sim t^{r},
\end{equation}
as $ t \to 0 $. Then we say that the fields $ X_{1} $, $ X_{2} $, $
X_{3} $ of (\ref{6.3.4}) are in case $ I_{r} $ or of type $ I_{r} $, $
r > 0 $. This implies that
$$ 
X_{3}(\gamma_{(\bar{x}', \bar{\xi}')}(t))
\Big|_{X_{2}(\gamma_{(\bar{x}', \bar{\xi}')} (t)) = 0} \sim t^{q-1+r},
$$
for $ t \to 0 $ and $ \bar{x}' $ in a small neighborhood of the origin.
\end{definition}
The first and second examples in Section 7.1 for Case I operators are
of type $ I_{0} $, while the third example is of type $ I_{q-p} $. 

We will find this property to be relevant for the Gevrey hypoellipticity 
threshold of the corresponding sums of squares operators.

\section{Gevrey regularity for sums of squares of vector fields of
  type $ I_{0} $}
\setcounter{equation}{0}
\setcounter{theorem}{0}
\setcounter{proposition}{0}  
\setcounter{lemma}{0}
\setcounter{corollary}{0} 
\setcounter{definition}{0}
\renewcommand{\thetheorem}{\thesection.\arabic{theorem}}
\renewcommand{\theproposition}{\thesection.\arabic{proposition}}
\renewcommand{\thelemma}{\thesection.\arabic{lemma}}
\renewcommand{\thedefinition}{\thesection.\arabic{definition}}
\renewcommand{\thecorollary}{\thesection.\arabic{corollary}}
\renewcommand{\theequation}{\thesection.\arabic{equation}}
%

In this
section our purpose is to deduce microlocal Gevrey estimates for operators
of type
$ I_{0} $. For the sake of simplicity we slightly modify our notation in
(\ref{6.3.4}). Thus let us consider three vector fields of the form
\begin{eqnarray}
\label{9.1}
X_{1}(x, D) & = & D_{1} \nonumber \\
X_{2}(x, D) & = & a_{21}(x) D_{1} + x_{1}^{p-1} f_{2}(x) D_{2} +
x_{1}^{q-1} g_{2}(x) D_{3}  \\
X_{3}(x, D) & = & a_{31}(x) D_{1} + x_{1}^{p-1} f_{3}(x) D_{2} +
x_{1}^{q-1} g_{3}(x) D_{3},
\nonumber
\end{eqnarray}
where $ f_{j} $ and $ g_{j} $, $ j = 2, 3 $, are real analytic
functions defined in a neighborhood of the origin and such that
(\ref{6.3.5:caseIlastlayer}) becomes
\begin{eqnarray}
\label{9.2}
f_{2}(0, x') & \neq & 0  \\ 
f_{3}(0, x') & = & \lambda(x') f_{2}(0, x') \nonumber 
\end{eqnarray}
and 
\begin{equation}
\label{9.3}
\begin{bmatrix}
  g_{2} (0)  \\
  g_{3}(0) 
\end{bmatrix}
\neq 0.
\end{equation}
Moreover (\ref{6.3.6:case1linind}) becomes
\begin{equation}
\label{9.4}
-\lambda(x') g_{2}(x) + g_{3}(x) + \frac{x_{1}}{f_{2}(0, x')} \det 
\begin{bmatrix}
  \frac{f_{2} - f_{2}(0, x')}{x_{1}}  & g_{2}(x) \\
  \frac{f_{3} - f_{3}(0, x')}{x_{1}} &  g_{3}(x)
\end{bmatrix}
\neq 0
\end{equation}
if $ x_{1} \neq 0 $. Now the assumption that our operator $ \sum_{j =
  1}^{3} X_{j}^{2} $ is of type $ I_{0} $ means that
\begin{equation}
\label{9.5}- \lambda(0) g_{2}(0) + g_{3}(0) \neq 0.
\end{equation}
The latter implies (\ref{9.4}), while (\ref{9.4}) makes sense due to
(\ref{9.3}).

\begin{lemma}
\label{lemma6.1}
Let $ \alpha $, $ \beta $ and $ \gamma $ be real analytic
functions defined in a neighborhood of the origin in $ \R^{3} $. Then we
can find real analytic functions $ a $, $ b $ and $ c $ such that
$$ 
\alpha(x) D_{1} + \beta(x) x_{1}^{p-1} D_{2} +  \gamma(x) x_{1}^{q-1}
D_{3} = a(x) X_{1} + b(x) X_{2} + c(x) X_{3}.
$$
\end{lemma}
\begin{proof}
This very useful lemma is a simple consequence of the assumptions, and says
that the span of the vector fields $\{D_1, x^{p-1}D_2, x_{q-1}D_3\}$ is that
same as that of the vector fields $\{X_j\}.$ Using elementary row and
column operations on the matrix on the right hand side of (\ref{9.1}) the
Lemma states the invertibility of the matrix 
\begin{equation}
\label{9.6}
\begin{bmatrix}
  f_{2}(x)  & g_{2}(x) \\
  f_{3}(x)  &  g_{3}(x)
\end{bmatrix}
\end{equation}
which, in view of (\ref{9.2}) is equivalent to the invertibility of the matrix
\begin{equation}
\label{9.7}
\begin{bmatrix}
  1  & g_{2}(x) \\
  0  &  g_{3}(x) - \lambda(x)g_2(x)
\end{bmatrix}.
\end{equation}
But this is just (\ref{9.5}) (all locally).
\end{proof}

\begin{lemma}
\label{lemma9.2}
For $j=1,2,3,$ and $m$ an integer, 
\begin{equation}
\label{9.8}
[X_{j} , D_{3}^{m} ] = \sum_{\ell =1}^{m} \binom{m}{\ell}
\sum_{h=1}^{3} \tilde{\gamma}_{jh}^{(\ell)} X_{h} D_{3}^{m-\ell},
\end{equation}
where 
\begin{equation} 
\label{9.9}
| \partial^{\alpha}   \tilde{\gamma}_{j h}^{(\ell)} |
\lesssim C_{j h}^{\ell + |\alpha|} (\ell + |\alpha|)!
\end{equation} 
Equivalently, 
\begin{equation}
\label{9.10}
D_{3}^{m} X_{j} = - \sum_{\ell =0}^{m} \binom{m}{\ell}
\sum_{h=1}^{3} \tilde{\gamma}_{jh}^{(\ell)} X_{h} D_{3}^{m-\ell},
\end{equation}
where $ \tilde\gamma_{jh}^{(0)} = - \delta_{jh} $. 
\end{lemma}
\begin{proof}
This is just an iteration of the previous Lemma. 
\end{proof}

Let denote by $ \phi $ a cut off function identically equal to
one in a neighborhood of the origin in $ \R^{3} $. Due to the
special form of our coordinates and the fact that the characteristic
manifold is simplectic, we may assume that $ \phi $ is independent of
the variable $ x_{1} $: in fact we may always take $ \phi $ as a product of
three such cut off functions each depending on a single coordinate, $
x_{j} $, and every $ x_{1} $-derivative landing on $ \phi(x_{1}) $
would leave a cut off supported in a region where $ x_{1}$ is bounded 
away from zero, hence in a regioln 
where the operator is (uniformly, microlocally) elliptic.
Thus we take $ \phi(x) = \phi(x')$. Here $
\phi$ is assumed to be a function of Ehrenpreis-H\"ormander type (see
e.g. \cite{Ehrenpreis1960}, \cite{H1971}), i.e.,
denoting by $ U $ our neighborhood of the origin, then $ \phi_{j} $
has the following property: for any $ \tilde{U}$ compactly contained
in $ U $, and for any fixed $ r \in \N $, we choose $ \phi_{j} = \phi_{j,
r}
\in C_{0}^{\infty}(U)$, $ \phi \equiv 1 $ on $ \tilde{U} $
and such that, with a {\em universal} constant (i.e.,
depending only
on the dimension of the Euclidean space in which we work)
$ C_{0} $ 
such that 
$$ |\phi^{(k)}(x) | \leq 
\left(\frac{C_{0}}{\text{dist}(\tilde{U},U^c)}\right)^{k+1} r^k 
\text{ for } k \leq 3r. 
$$
It is a well known fact that the operator
%
%
$$ 
P (x, D) = \sum_{j=1}^{3} X_{j}(x, D)^{2}
$$
is $ C^{\infty} $-hypoelliptic and satisfies an a priori estimate of
the form
\begin{equation}
\label{9.11}
\sum _{j=1}^{3} \| X_{j} u \|^{2} + \| u \|_{1/q}^{2} \leq C \left(
  |\langle P u, u \rangle| + \| u \|^{2} \right),
\end{equation}
where $ u $ is a rapidly decreasing smooth function, $ \| \cdot \|_{s}
$ denotes the usual Sobolev norm of order $ s $ and $ \| \cdot \| = \|
\cdot \|_{0}$ is the $ L^{2} $ norm.

We want to obtain a bound for an expression of the form
\begin{equation}
\label{9.12}
\| X_{j} \phi(x') D_{3}^{r} u \|,
\end{equation}
where, since we are in a microlocal neighborhood of the point $ (0,
e_{3}) $, $ D_{3} $ is an elliptic operator. It is well known that
obtaining a bound for (\ref{9.12}) of the type $ \| X_{j} \phi(x')
D_{3}^{r} u \| \leq C^{r+1} r!^{s} $ allows us to deduce that $ P $ is
Gevrey (micro-)hypoelliptic of order $ s $.

\begin{remark}We would like to mention here that in the case of
the second example of Section 7.1, i.e. the
Oleinik-Radkevi\v c\ operator, the authors in
\cite{BT} proved that one has $G^{q/p}$ hypoellipticity and that this
bound is optimal.
\end{remark}

Instead of bounding the quantity in (\ref{9.12}), for technical
reasons we want to bound the more general quantity:
\begin{equation}
\label{9.13}
\| X_{j} x_{1}^{a} \phi^{(b)} D_{3}^{r-c} u \| + \| x_{1}^{a}
\phi^{(b)} D_{3}^{r-c} u \|
_{1/q}\end{equation}
where $ a $, $ b $ and $ c $ are positive integers with $a\leq q$ but 
$ b $ and $ c $
bounded only by $r.$ 
Using (\ref{9.11}), we see that (\ref{9.13}) is
bounded by 
\begin{equation}
\label{9.14}
| \langle P x_{1}^{a} \phi^{(b)}
D_{3}^{r-c} u,\; x_{1}^{a} \phi^{(b)} D_{3}^{r-c} u \rangle | + \|
x_{1}^{a} \phi^{(b)} D_{3}^{r-c} u \|^{2},  
\end{equation}
modulo a positive constant in front of everything appearing in the
above formula. We need to move $ P $ in (\ref{9.14}) to the right
(onto $ u $); the term with the $ L^{2} $ norm will be easier to
handle. Writing $P=\Sigma X^2 $ and then $[X^2,V]=X[X,V]+[X,V]X$
with $V=x_{1}^{a} \phi^{(b)} D_{3}^{r-c},$ we find

$$
\langle P x_{1}^{a} \phi^{(b)} D_{3}^{r-c} u , x_{1}^{a}
\phi^{(b)} D_{3}^{r-c} u \rangle =   
\langle x_{1}^{a} \phi^{(b)}
D_{3}^{r-c} P u , x_{1}^{a} \phi^{(b)} D_{3}^{r-c} u
\rangle $$
\begin{multline}
\label{9.15}
\qquad \qquad \qquad + \sum_{j=1}^3\langle X_j[ X_{j}, x_{1}^{a}
\phi^{(b)} D_{3}^{r-c}] u ,
 x_{1}^{a} \phi^{(b)} D_{3}^{r-c} u \rangle \\
 + \sum_{j=1}^3\langle [ X_{j}, x_{1}^{a}
\phi^{(b)} D_{3}^{r-c}] X_{j} u ,  x_{1}^{a} \phi^{(b)} D_{3}^{r-c} u
\rangle.
\end{multline}
The first of the above right hand side terms is good, since we assume
$ Pu $ to be analytic, even $0$. 

The second and third terms on the
right hand side in (\ref{9.15}) have many common features, which we
may treat with the help of Lemma \ref{lemma9.2}. 

For j=1, we have
\begin{equation}
\label{9.16}
[ X_{1}, x_{1}^{a} \phi^{(b)} D_{3}^{r-c}] =
ax_1^{a-1}\phi^{(b)}D_3^{r-c},
\end{equation}
and so
\begin{multline}
\label{9.17}
[ X_{1}, x_{1}^{a} \phi^{(b)} D_{3}^{r-c}]X_1 =
X_1ax_1^{a-1}\phi^{(b)}D_3^{r-c}  
- a(a-1)x_1^{a-2}\phi^{(b)}D_3^{r-c}.
\end{multline}
For $j=2,3,$
\begin{multline}
\label{9.18}
[ X_{j}, x_{1}^{a} \phi^{(b)} D_{3}^{r-c}] =
f(x)x_1^{p-1}x_1^{a}\phi^{(b+1)}D_3^{r-c} + \\ + x_{1}^{a} \phi^{(b)}
\sum_{\ell =1}^{r-c} \binom{r-c}{\ell}
\sum_{h=1}^{3} \tilde{\gamma}_{jh}^{(\ell)} X_{h} D_{3}^{r-c-\ell},
\end{multline}
with $f$ analytic, and thus using Lemma \ref{lemma9.2} again, ($j=2,3$)
\begin{multline}
\label{9.19}
[ X_{j}, x_{1}^{a} \phi^{(b)} D_{3}^{r-c}]X_j =
f(x)x_1^{p-1}x_1^{a}\phi^{(b+1)}D_3^{r-c}X_j + \\ + x_{1}^{a}
\phi^{(b)}
\sum_{\ell =1}^{r-c} \binom{r-c}{\ell}
\sum_{h=1}^{3} \tilde{\gamma}_{jh}^{(\ell)} X_{h} D_{3}^{r-c-\ell}X_j =
\\ = f(x)x_1^{p-1}x_1^{a}\phi^{(b+1)}\sum_{\ell =0}^{r-c}
\binom{r-c}{\ell}
\sum_{h=1}^{3} \tilde{\gamma}_{jh}^{(\ell)} X_{h} D_{3}^{r-c-\ell} +
\\ + x_{1}^{a}
\phi^{(b)}
\sum_{\ell =1}^{r-c} \binom{r-c}{\ell}
\sum_{h=1}^{3} \tilde{\gamma}_{jh}^{(\ell)} X_{h} \sum_{\ell_1
=0}^{r-c-\ell}
\binom{r-c-\ell}{\ell_1}
\sum_{k=1}^{3} \tilde{\gamma}_{jk}^{(\ell_1)} X_{k}
D_{3}^{r-c-\ell-\ell_1}.
\end{multline}
Going back to (\ref{9.15}), the first term we have seen is harmless as it
contains $Pu.$ In the second, we integrate by parts and use a weighted
Schwarz inequality. Since $X_j^*$ is equal to $-X_j$ modulo a zero
order term,  
the second term on the right in (\ref{9.15}), using (\ref{9.16}) and 
(\ref{9.18}), becomes  
\begin{multline}
\label{9.20}
|\sum_{j=1}^3\langle X_j[ X_{j}, x_{1}^{a} \phi^{(b)}
D_{3}^{r-c}] u ,
 x_{1}^{a} \phi^{(b)} D_{3}^{r-c} u \rangle| \leq \\
\leq \epsilon
\sum_{k=1}^3\| X_kx_{1}^{a}
\phi^{(b)} D_{3}^{r-c} u\|^2 + C_\epsilon\| a
x_1^{a-1}\phi^{(b)}D_3^{r-c}u\|^2 + \\ + C_\epsilon
\|f(x)x_1^{p-1}x_1^{a}\phi^{(b+1)}D_3^{r-c}u\|^2 +\\+ C_\epsilon
\sum_{1\leq h\leq 3 \atop 2\leq j\leq 3}\| x_{1}^{a} \phi^{(b)}
\sum_{\ell =1}^{r-c} \binom{r-c}{\ell}
\tilde{\gamma}_{jh}^{(\ell)} X_{h}
D_{3}^{r-c-\ell}u\|^2. 
\end{multline}
This expression we leave for the moment and treat the issues which arise
in the double commutator needed for the last term in (\ref{9.15}), those
which have already been expanded in (\ref{9.17}) and (\ref{9.19}).

We may continue with (\ref{9.17}) in (\ref{9.15}):
\begin{multline}
\label{9.21}
\langle [ X_{1}, x_{1}^{a} \phi^{(b)} D_{3}^{r-c}] X_{1} u , 
x_{1}^{a} \phi^{(b)} D_{3}^{r-c} u \rangle = \\
=
\langle X_1ax_1^{a-1}\phi^{(b)}D_3^{r-c} u , 
x_{1}^{a} \phi^{(b)} D_{3}^{r-c} u \rangle \\
-
\langle a(a-1)x_1^{a-2}\phi^{(b)}D_3^{r-c} u , 
x_{1}^{a} \phi^{(b)} D_{3}^{r-c} u \rangle .  
\end{multline}
We shall also continue with (\ref{9.19}) in (\ref{9.15}): for $j=2,3$
\begin{multline}
\label{9.22}
\langle [ X_{j}, x_{1}^{a} \phi^{(b)} D_{3}^{r-c}] X_{j} u , 
x_{1}^{a} \phi^{(b)} D_{3}^{r-c} u \rangle = \\
=
\sum_{\ell =0\atop h=1,2,3}^{r-c}
\binom{r-c}{\ell}
\langle f(x)\tilde{\gamma}_{jh}^{(\ell)}x_1^{p-1}x_1^{a}\phi^{(b+1)} 
X_{h} D_{3}^{r-c-\ell} u ,  x_{1}^{a} \phi^{(b)} D_{3}^{r-c} u \rangle \\
 + \sum_{\ell +\ell_1=1}^{r-c}
\binom{r-c}{\ell, \ell_1}\langle x_{1}^{a}
\phi^{(b)}
\sum_{h,k=1}^{3} \tilde{\gamma}_{jh}^{(\ell)} X_{h} 
 \tilde{\gamma}_{jk}^{(\ell_1)} X_{k}
D_{3}^{r-c-\ell-\ell_1} u , 
x_{1}^{a} \phi^{(b)} D_{3}^{r-c} u \rangle .
\end{multline}
Here we have used the `multinomial' notation for brevity:
$${\alpha \choose \beta, \gamma} = \frac{\alpha !}{\beta !
\gamma ! (\alpha-\beta-\gamma)!}={\alpha \choose
\beta}{\alpha-\beta \choose \gamma}.$$

Before collecting our individual terms we throw in a kind of
`symmetrization' of the first term on the left, for errors will often
appear in this form. In so doing, we will encounter one more
commutator, which is covered under the fourth and fifth terms on
the right, hence contributing nothing new. We also drop the
subscripts on the vector fields now. From (\ref{9.13}), (\ref{9.14}), and
(\ref{9.15}), (\ref{9.20}), (\ref{9.21}), (\ref{9.22}) we have, for any
positive $\epsilon$,
\begin{multline}
\label{9.23}
\| X x_{1}^{a} \phi^{(b)} D_{3}^{r-c} u \| + 
\| x_{1}^{a} \phi^{(b)} D_{3}^{r-c} u \|_{1/q} +  \|  x_{1}^{a}
\phi^{(b)} XD_{3}^{r-c} u \| \lesssim \\
\lesssim
\|x_{1}^{a} \phi^{(b)} D_{3}^{r-c} P u\| + \|
x_{1}^{a} \phi^{(b)} D_{3}^{r-c} u \|  \\
+\epsilon \| X x_{1}^{a}
\phi^{(b)} D_{3}^{r-c} u\| + C_\epsilon\| a
x_1^{a-1}\phi^{(b)}D_3^{r-c}u\| + \\ + C_\epsilon
\|f(x)x_1^{p-1}x_1^{a}\phi^{(b+1)}D_3^{r-c}u\| +\\
+ C_\epsilon
\sum_{2\leq j\leq 3}\| x_{1}^{a} \phi^{(b)}
\sum_{\ell =1}^{r-c} \binom{r-c}{\ell}
\tilde{\gamma}_{jh}^{(\ell)} X
D_{3}^{r-c-\ell}u\| + \\
+|\langle Xax_1^{a-1}\phi^{(b)}D_3^{r-c} u , 
x_{1}^{a} \phi^{(b)} D_{3}^{r-c} u \rangle |^{1/2} \\
+|\langle a(a-1)x_1^{a-2}\phi^{(b)}D_3^{r-c} u , 
x_{1}^{a} \phi^{(b)} D_{3}^{r-c} u \rangle  |^{1/2}+ \\
+ \sum_{\ell =0\atop {j=2,3}}^{r-c}
\binom{r-c}{\ell}
|\langle f(x)\tilde{\gamma}_{jh}^{(\ell)}x_1^{p-1}x_1^{a}\phi^{(b+1)} 
XD_{3}^{r-c-\ell} u ,  x_{1}^{a} \phi^{(b)} D_{3}^{r-c} u \rangle
|^{1/2} \\
 + \sum_{\ell +\ell_1=1\atop j=2,3}^{r-c}
\binom{r-c}{\ell, \ell_1}|\langle x_{1}^{a}
\phi^{(b)}
\tilde{\gamma}_{j\cdot}^{(\ell)} X 
 \tilde{\gamma}_{j\cdot}^{(\ell_1)} XD_{3}^{r-c-\ell-\ell_1} u , 
x_{1}^{a} \phi^{(b)} D_{3}^{r-c} u \rangle|^{1/2} = \\
=
I_1+I_2+I_3+I_4+I_5+I_6+I_7+I_8+I_9.
\end{multline}
\subsection{The term $I_1$} This term is harmless since $Pu$ is real
analytic, even zero, in the support of all $\phi.$
\subsection{The term $I_2$} This term will be bounded by a small multiple
of (\ref{9.13}) if we take the support of all the
localizing functions small, and hence may be absorbed. 
\subsection{The term $I_3$} This term is already a small multiple of
(\ref{9.13}), hence absorbable for $\epsilon$ small.
\subsection{The term $I_4$} This term, $C_\epsilon\| a
x_1^{a-1}\phi^{(b)}D_3^{r-c}u\|,$ exhibits an overall gain (in the
norm) of $1/q,$ but pays for it with a decrease in the power of $x_1.$
We will consider this term further below. 
\subsection{The term $I_5$} This term, bounded at once by
$$C_fC_\epsilon
\|x_1^{p-1}x_1^{a}\phi^{(b+1)}D_3^{r-c}u\|,$$ 
suffers a new derivative on $\phi$ but gains the factor $x_1^{p-1},$ 
\subsection{The term $I_6$} This term, easily bounded by
$$\tilde{C}_\epsilon
\sup_{1\leq \ell \leq r-c}
(C_\gamma r)^\ell\| x_{1}^{a} \phi^{(b)}
 X D_{3}^{r-c-\ell}u\|$$
in view of the estimates (\ref{9.9}), where
$C_\gamma$ depends only on the coefficients of the $X_j$ and their first
few derivatives. This term will be further treated under $I_8$ below,
where also the term with $\ell = 0$ appears, though with a small constant
in front. 
\subsection{The term $I_7$} This term, 
$$|\langle a(a-1)x_1^{a-2}\phi^{(b)}D_3^{r-c} u , 
x_{1}^{a} \phi^{(b)} D_{3}^{r-c} u \rangle  |^{1/2},$$
is bounded exactly as is $I_4$ above once one power of $a$ is moved to
the left and the Schwarz inequality applied. 
\subsection{The term $I_8$} This term, 
$$\sum_{\ell =0\atop {j=2,3}}^{r-c}
\binom{r-c}{\ell}
|\langle
f(x)\tilde{\gamma}_{j\cdot}^{(\ell)}x_1^{p-1}x_1^{a}\phi^{(b+1)} 
X D_{3}^{r-c-\ell} u ,  x_{1}^{a} \phi^{(b)} D_{3}^{r-c} u 
\rangle |^{1/2},$$
permits us to move $x_1^{p-1}\phi^{(b+1)}$ to the right and
$\phi^{(b)}$ to the left, apply the Schwarz inequality and bring both $f$
and $\tilde{\gamma}_{j\cdot}^{(\ell)}$ out of the norm. The result is 
\begin{multline}
\epsilon\|x^a\phi^{(b)}XD_3^{r-c}u\| + \sup_{\ell \geq
1}(C_\gamma r)^\ell\|x^a\phi^{(b)}XD_3^{r-c-\ell}u\| 
\\ + C_\epsilon\|x_1^{a+p-1}\phi^{(b+1)}D_3^{r-c}u\|.
\end{multline}
The last of these is exactly like $I_5$ above, while the supremum, 
 has been met in $I_6$ above. 
The first term, which we note carries the small constant $\epsilon,$  will 
be absorbed on the left hand side of (\ref{9.13}) once the $X$ is
commuted to the left. 
\subsection{The term $I_9$} This term, 
$$\sum_{\ell +\ell_1=1\atop
j=2,3}^{r-c}
\binom{r-c}{\ell, \ell_1}|\langle x_{1}^{a}
\phi^{(b)}
\tilde{\gamma}_{j\cdot}^{(\ell)} X
 \tilde{\gamma}_{j\cdot}^{(\ell_1)} X
D_{3}^{r-c-\ell-\ell_1} u , 
x_{1}^{a} \phi^{(b)} D_{3}^{r-c} u \rangle|^{1/2},$$
carries with it some of the features of all of the above terms. We want to
move $X$ to the right, use the weighted Schwarz inequality, and
estimate the derivatives $\tilde{\gamma}_{j\cdot}^{(\ell)}$ just as we
have done before. But two things may happen: in first commuting $X$ to
the left another derivative may fall on $\tilde{\gamma}_{j\cdot}^{(\ell)}$
(doing no harm - the estimates on these derivatives are flexible enough to
handlle one or two more derivatives by changing the constant a bit,
uniformly in
$r$). But the coefficient $x_1^a\phi^{(b)}$ may also be differentiated by
$X.$ No matter - this has happened often before, and either $\phi$
receives one more derivative gains a coefficient of $x_1^{p-1},$
as in $I_5,$ or $x_1^a$ becomes $x_1^{a-1}$ as in $I_4.$

Putting these results together, the error terms, apart from those which
may be absorbed on the left, we have arrived at
\begin{lemma} 
\label{lemma9.3}
For any $a, b, c,$ and $r$ we have the
estimate
$$
\| X x_{1}^{a} \phi^{(b)} D_{3}^{r-c} u \| + 
\| x_{1}^{a} \phi^{(b)} D_{3}^{r-c} u \|_{1/q} +  \|  x_{1}^{a}
\phi^{(b)} XD_{3}^{r-c} u \| \\
\lesssim $$
\begin{multline}
\label{9.24}
\|x_{1}^{a} \phi^{(b)} D_{3}^{r-c} P u\| 
+ \sup_{1\leq \ell \leq r-c}
(C_\gamma r)^\ell\| x_{1}^{a} \phi^{(b)}
 X D_{3}^{r-c-\ell}u\| \\
+ \| x_1^{a-1}\phi^{(b)}D_3^{r-c}u\| + 
\|x_1^{a+p-1}\phi^{(b+1)}D_3^{r-c}u\| = J_1 + J_2 + J_3
+J_4.
\end{multline}
\end{lemma}
And these terms are of four distinct types: the first involves $Pu$
and is harmless; the second exhibits a gain of $\ell$ powers of $D_3$
at the expense of $\ell$ powers of $r;$ iteration will lead to $(Cr)^r
\sim \tilde{C}^r r!,$ which by itself would lead to analytic growth.  

For the final two terms, $J_3$ and $J_4,$ we argue as follows: 
\begin{itemize}
\item[1)] In treating terms where a power of $x$ has been
differentiated, we invoke subellipticity, writing
$$
\| x_1^{a-1}\phi^{(b)}D_3^{r-c}u\| = 
\|x_1^{a-1}\phi^{(b)}D_3^{r-c-1/q}u\|_{1/q} + E
$$
and estimate $E$ using the standard calculus of
pseudo-differential operators - giving rise to a sum of terms,
in which a typical term has $k$ more derivatives on $\phi$
and $k$ fewer powers of $D_3,$ modulo an error with no
derivatives on $u.$ This trade-off, $D_3$'s being
transferred from $u$ to $\phi,$ is the sort that would lead
to analyticity. At any rate, the principal contribution is
similar to the second term on the left of
Lemma \ref{lemma9.3} with
$a$ decreased by one and $c$ increased by $1/q.$
\item[2)] We observe that when $a=0$ (at the outset, for
instance) this kind of term does not arise; thus $J_4$ will be
the first term to arise, starting from $\|\phi
D_3^{r}u\|:
\|\phi D_3^{r}u\| \rightarrow
\|x_1^{p-1}\phi^{(1)} D_3^{r-1/q}u\|_{1/q}.$ 

\item[3)] Alternatively, when, as in 2) just above,  one does
add $p-1$ powers of
$x_1,$ add a full derivative to $\phi,$ one {\em may}
reach a total of $q-1\ \ x$'s, in which
case one invokes Lemma \ref{lemma6.1}, 
and writes $x^{q-1}D_3 = \sum
b_jX_j$ and does {\em not} employ the $1/q-$ `shunt' in the first item
just above. If the powers of $x$ do not permit this, we use 
the subellipticity again.  

\item[4)] All together, then, we observe that after $s$ steps
of type 1) and
$k$ steps of type 3), starting from $a=b=c=0,$ we will have, as
`worst' errors, 
$$
C^{k+s}r^k\| x_{1}^{k(p-1) - s} \phi^{(k)} D_{3}^{r-\ell-
\frac{k+s-1}{q}} u \|,
$$
where after the last step we have {\em not} taken $D_3^{1/q}$ and
moved it to be part of the norm; for this time, assuming that we have
approximately $q-1$ powers of $x_1,$  we will use Lemma
\ref{lemma6.1} to `create' an $X.$
\end{itemize}

Whenever possible (when the powers of $x_1$ grow to $q-1,$ we do not
take advantage of the subelliptic $1/q$ gain but combine $x^{q-1}$ with
$D_3$ to produce an $X$ instead. This may happen $t$ times. The result
is that after $s+k+t$ iterations we have an expression  
\begin{equation}
\label{9.25}
r^{\ell} \| \phi^{(k)} x_{1}^{k(p-1) -s -t (q-1)} D_{3}^{r -\ell -
  \frac{k+s-t}{q} - t} u \|.
\end{equation}
Now $| \phi^{(k)}| \leq C^{k+1} r^{k} $, so that (\ref{9.25})
is bounded by
\begin{equation}
\label{9.26}
C^{k+1} r^{\ell+k} \|x_{1}^{k(p-1) -s -t (q-1)} D_{3}^{r -\ell -
  \frac{k+s-t}{q} - t} u \|_{L^2(supp \phi)}.
\end{equation}
Since we are looking for powers of $ x_{1} $ as close to zero as
possible (where we started) to gauge the effect of returning to the
starting point, it is natural to take 
$$ 
t = \frac{k(p-1) -s}{q-1},
$$
or its integer part. This choice of $ t $ reduces the quantity in
(\ref{9.26}) to the following
\begin{equation}
\label{9.27}
C^{k+1} r^{\ell+k} \| D_{3}^{r - \left( \ell + k \frac{p}{q}\right)} u
\|. 
\end{equation}
Upon iteration we get
\begin{equation}
\label{9.28}
C^{\sum_{j} k_{j}} C^{r} r^{\sum_{j} \ell_{j}+ \sum_{j} k_{j} } \|
D_{3}^{r - \sum_{j}\ell_{j} - \frac{p}{q} \sum_{j} k_{j}} u \|, 
\end{equation}
where
$$ 
r - \sum_{j}\ell_{j} - \frac{p}{q} \sum_{j} k_{j} \sim 0.
$$
Let us write $ K = \sum_{j} k_{j} $ and $ L = \sum_{j} \ell_{j}
$. Then
$$ 
K + L = r \frac{K+L}{r} \sim r \frac{ K+L}{ L +
  \frac{\textstyle p}{\textstyle q} K} \leq r 
\frac{q}{p} , 
$$
since $ q \geq p $. This ends the proof of the following 
\begin{theorem}
\label{theorem2}
The operator 
$$ 
P (x, D) = \sum_{j=1}^{3} X_{j}^{2},
$$
where the $ X_{j} $ are given in (\ref{9.1}), is Gevrey hypoelliptic
of order $ \frac{q}{p} $, i.e.
$$ 
P u = f \in C^{\omega}\ \  \text{implies that microlocally}\ \  u \in G^{s},
\quad s \geq \frac{q}{p}.
$$
\end{theorem}
\noindent {\bf Remarks}
The results given are microlocal. To provide a proof in all detail would
entail introducing cut-off functions which are local in space $x$ and
also in the frequency variables $\xi$.  This can be done, and has been
carried out in all detail in
\cite{Tartakoff1980} and \cite{Tartakoff1981} in the analytic case
and in \cite{BTprop} in the Gevrey category. One introduces
localizing functions with the local behavior used here and conic
localization in the frequency variables, all cut-off near the origin
in the dual variables (with analytic error) in the manner detailed in
\cite{Tartakoff1981}. We omit details here, as they would largely
repeat \cite{Tartakoff1981} and risk rendering the exposition
unreadable. 

In the case of the Oleinik-Radkevich
model we know from \cite{BT} that these results are optimal and that in 
particular the result is analytic hypoelliptic if and only if $p=q.$ We
also strongly believe, but have not yet been able to prove, that every
`threshold' obtained in this paper is also sharp.

\end{document}